\documentclass[11pt, a4paper, titlepage]{article}

\usepackage{Stylefile}

\usepackage[T1]{fontenc}
\usepackage[utf8]{inputenc}
\usepackage{lmodern}
\usepackage{textcomp}
\usepackage{microtype}

\usepackage[english]{babel}
\usepackage{csquotes}

\usepackage[backend=biber,
			style=ext-numeric,
			isbn=false,
			doi=false,
            maxbibnames=99,
            sortcites,
            articlein=false
			]{biblatex}
\usepackage{graphicx}
\usepackage{setspace}
\usepackage{caption}
\usepackage{subcaption}

\usepackage{amsmath}
\usepackage{amsfonts}
\usepackage{amssymb}
\usepackage{amsthm}
\usepackage{multicol}
\usepackage{tikz}
\usetikzlibrary{arrows.meta, decorations.markings}
\usepackage{pgfplots}
\usepackage{pdfpages}
\pgfplotsset{compat=newest}
\usepackage{nicematrix}
\usepackage{mathtools}
\usepackage{xcolor}
\usepackage{hhline}

\numberwithin{equation}{section}

\newcommand*{\Z}{\mathbb{Z}}

\newcommand*{\R}{\mathbb{R}}

\newcommand*{\eps}{\varepsilon}
\newcommand{\comm}[1]{\textcolor{black}{#1}}

\usepackage{thmtools}

\declaretheorem[
	name=Theorem,
	numberwithin=section
	]{thm}
\declaretheorem[
	name=Lemma,
	sibling=thm,
	]{lem}

\declaretheorem[
	name=Definition,
	sibling=thm,
	]{defin}

\declaretheorem[
	name=Remark,
	sibling=thm
	]{rem}

\usepackage[pdftex,pdfpagelabels]{hyperref}

\begin{document}


\makeatletter
\tikzset{
    dot diameter/.store in=\dot@diameter,
    dot diameter=2pt,
    dot spacing/.store in=\dot@spacing,
    dot spacing=8pt,
    dots/.style={
        line width=\dot@diameter,
        line cap=round,
        dash pattern=on 0pt off \dot@spacing
    }}
\tikzset{->-/.style={decoration={
  markings,
  mark=at position #1 with {\arrow{{Straight Barb}}}},postaction={decorate}}}

\tikzset{->>-/.style={decoration={
  markings,
  mark=at position #1 with {\arrow{{Straight Barb} {Straight Barb}}}},postaction={decorate}}}

\makeatother

\renewcommand{\shorttitle}{The Robertson Model}

\hypersetup{
pdftitle={A Multi-Parameter Singular Perturbation Analysis of the Robertson Model},
pdfauthor={L.~Baumgartner, P.~Szmolyan},
pdfkeywords={Multi-parameter singular perturbations, Robertson model, Geometric singular perturbation theory, Blow-up method},
}

\vspace*{1.2cm}

\begin{center}
\LARGE A Multi-Parameter Singular Perturbation Analysis\\ of the Robertson Model
\end{center}

\begin{center}
    \large L. Baumgartner\footnote{Institute of Analysis and Scientific Computing, TU Wien, Wiedner Hauptstraße 8–10, 1040 Vienna, Austria. \\
    E-Mail: lukas.baumgartner@tuwien.ac.at} \& P. Szmolyan\footnote{Institute of Analysis and Scientific Computing, TU Wien, Wiedner Hauptstraße 8–10, 1040 Vienna, Austria. \\
    E-Mail: peter.szmolyan@tuwien.ac.at} 
\end{center}

\onehalfspacing 

\begin{abstract}
	The Robertson model describing a chemical reaction involving three reactants is one of the classical examples of stiffness in ODEs. The stiffness is caused by the occurrence of three reaction rates $k_1,\,k_2$, and $k_3$, with largely differing orders of magnitude, acting as parameters. The model has been widely used as a numerical test problem. Surprisingly, no asymptotic analysis of this multiscale problem seems to exist. In this paper we provide a full asymptotic analysis of the Robertson model under the assumption $k_1, k_3 \ll k_2$. We rewrite the equations as a two-parameter singular perturbation problem in the rescaled small parameters $(\eps_1,\eps_2):=(k_1/k_2,k_3/k_2)$, which we then analyze using geometric singular perturbation theory (GSPT). To deal with the multi-parameter singular structure, we perform blow-ups in parameter- and variable space. We identify four distinct regimes in a neighbourhood of the singular limit \mbox{$(\eps_1,\eps_2)= (0,0)$}. Within these four regimes we use GSPT and additional blow-ups to analyze the dynamics and the structure of solutions. Our asymptotic results are in excellent qualitative and quantitative agreement with the numerics.

\end{abstract}

\keywords{Multi-parameter singular perturbation \and Robertson model \and Geometric singular perturbation theory \and Blow-up method}
\vspace{-0.3cm}
\classifications{34E10 \and 34E13 \and 34E15 \and 92E20}

\section{Introduction} \label{sec:Introduction}

\begin{figure}[h]
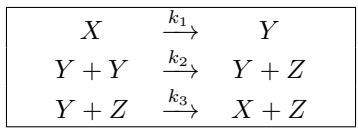

\centering
\begin{tabular}{|c c c c c|}
\hline 
 & $X$ & $\xrightarrow{K_1}$ & $Y$ &\\
 & $Y+Y$ & $\xrightarrow{K_2}$ & $Y+Z$ &\\
 & $Y+Z$ & $\xrightarrow{K_3}$ & $X+Z$ &\\
\hline
\end{tabular}
\caption{Reaction scheme of the Robertson model.}
\label{fig:rob_reaction_scheme}
\end{figure}

In this paper we give a dynamical systems analysis of the Robertson model \cite{robertson_solution_1966} based on methods from geometric singular perturbation theory (GSPT). The Robertson model describes a chemical reaction  of three reactants $X,\,Y$, and $Z$, which interact according to the reaction scheme shown in Figure \ref{fig:rob_reaction_scheme}.

With mass-action kinetics the Robertson model leads to the following system of ODEs 

\begin{equation} \label{equ:rob_reaction}
\begin{aligned}
\dot{x}&=-k_1 x +k_3 yz\\
\dot{y}&=k_1 x - k_2 y^2 - k_3 yz\\
\dot{z}&=k_2 y^2,
\end{aligned}
\end{equation}

with corresponding concentrations $x,y,z \in \R$, reaction rates $k_i >0$, $i=1,2,3$. As usual\comm{,} $\dot{(\quad)}:=\frac{d}{dt}$ denotes the time derivative. The classical choice of parameters and initial values in \cite{robertson_solution_1966} is 
\begin{equation} \label{equ:rob_rates}
k_1=4\cdot10^{-2}, \quad k_2=3\cdot10^{7}, \quad k_3=10^4
\end{equation}
and 
\begin{equation} \label{eq:rob_initial}
    (x(0),y(0),z(0))^T=(1,0,0)^T. 
\end{equation}

The qualitative dynamics of system  \eqref{equ:rob_reaction} is fairly simple.

\begin{lem} \label{lem:robertson_textbook}
    All solutions of \eqref{equ:rob_reaction} starting in the non-negative orthant $\R^3_+$ exist globally in forward time. The $z$-axis is a line of attracting equilibria. The solution with initial value $(x_0,y_0,z_0)^T \in \R^3_+$ converges to the equilibrium $(\hat{x},\hat{y},\hat{z})^T=(0,0,c)^T$, with  
    $ c:=x_0+y_0+z_0 > 0$.

\end{lem}
\begin{proof}
    Adding the three equations of \eqref{equ:rob_reaction} implies that the quantity $x+y+z=const.$ is conserved. Since on the boundary of the non-negative orthant $\R^3_+$, the flow does not point outwards, i.e., 
$$\Dot{x}|_{x=0}=k_3yz\geq 0, \quad \Dot{y}|_{y=0}=k_1x\geq 0, \quad \Dot{z}|_{z=0}=k_2y^2\geq 0,$$
we can conclude that $\R^3_+$ is forward invariant under \eqref{equ:rob_reaction}, see \cite[p.~219]{Amann_1990}. Consequently, the solution starting at an initial value $(x_0,y_0,z_0)^T \in \R^3_+$, where $0 < c:=x_0+y_0+z_0$, is contained in the compact set 
$$K=\{(x,y,z)^T \in \R^3_+:x+y+z=c\}$$
and therefore exists for all times $t\geq 0$. 

Due to the conserved quantity, we may reduce the dimension of \eqref{equ:rob_reaction} by using $x=c-y-z$ to obtain
\begin{equation} \label{eq:rob_2D_intro}
    \begin{aligned}
        \dot{y}&=k_1 (c-y-z) - k_2 y^2 - k_3 yz\\
        \dot{z}&=k_2 y^2.
    \end{aligned}
\end{equation}
Since the divergence of the vector field \eqref{eq:rob_2D_intro} given by
$-k_1-2k_2y-k_3z$
is negative for positive reaction rates, we can exclude non-constant periodic solutions by the Bendixson-Dulac criterion. The unique equilibrium of \eqref{eq:rob_2D_intro} is given by $(\hat{y},\hat{z})^T=(0,c)^T$, hence by the Poincare-Bendixson theorem all solutions of \eqref{eq:rob_2D_intro} will ultimately converge to this equilibrium. 
\end{proof}

\begin{figure}
\centering
\includegraphics[scale=0.5]{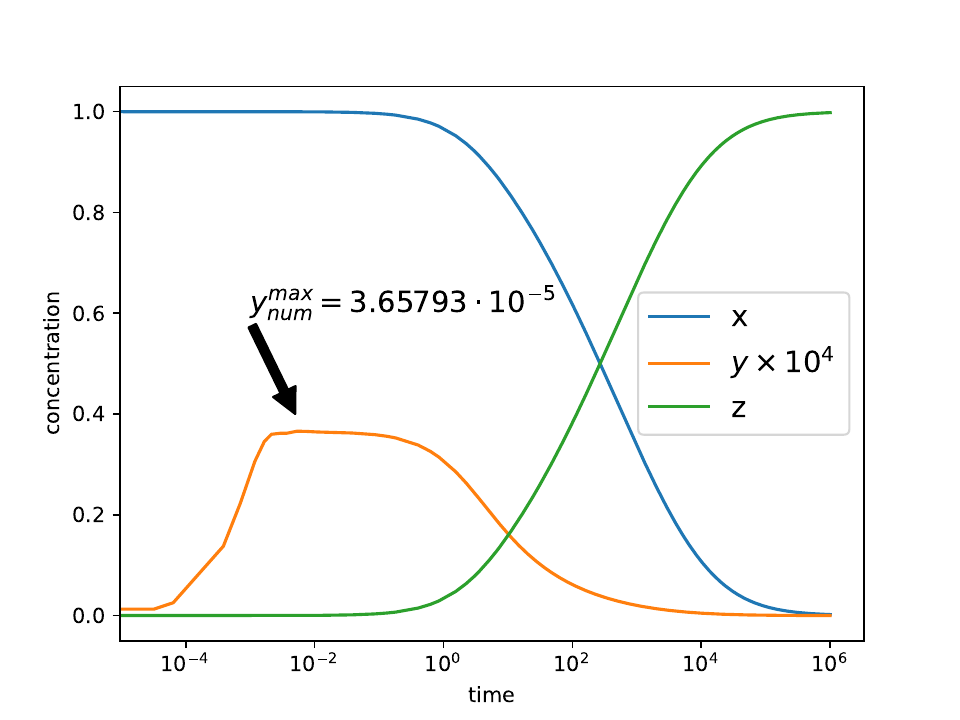}
\caption{Numerical solution of equation \eqref{equ:rob_reaction} with an implicit BDF solver \cite{Scipy_BDF}. Note the logarithmic time scale.}
\label{fig:rob_numeric}
\end{figure}

In particular, we conclude from Lemma \ref{lem:robertson_textbook} that the solution of \eqref{equ:rob_reaction} with initial value \eqref{eq:rob_initial} converges to the unique equilibrium $(\hat{x},\hat{y},\hat{z})^T=(0,0,1)^T$. Thus, our interest in the Robertson model is not this rather simple dynamics but the multi-scale structure of these solutions which we now describe in a preliminary way based on numerical simulations.
The time series of a numerical solution of \eqref{equ:rob_reaction} with the classical choice of reaction rates \eqref{equ:rob_rates} and initial condition \eqref{eq:rob_initial} is shown in Figure \ref{fig:rob_numeric}. 

In the time series\comm{,} three distinct parts can be distinguished. The reaction starts with a very fast initial increase of $y$ up to a plateau value $y^{max}_{num} \approx 3.65 \cdot 10^{-5}$. This is followed by an intermediate phase where $y$ is almost constant. In the third part the conversion of $x$ into $z$ (via $y$) proceeds on a much longer time scale. 

\begin{rem} \label{rem:early_comparison_y_max}
    \comm{Our analysis predicts a plateau value of $y^{max}=\sqrt{\frac{k_1 c}{k_2}}+\mathcal{O}(\frac{k_1}{k_2})$. A detailed derivation of this value can be found at the end of Section \ref{sec:Regime R2}. With the parameter values \eqref{equ:rob_rates} and initial condition \eqref{eq:rob_initial} (which implies $c=1$) this leads to $y^{max}=\frac{2}{\sqrt{3}}\cdot 10^{-\frac{9}{2}}+\mathcal{O}(10^{-9})=3.651\cdot 10^{-5}+\mathcal{O}(10^{-9})$, which matches the numerical value $y^{max}_{num}$ above.} 
\end{rem}

Numerical experiments indicate that this solution structure occurs for all parameter values 
\begin{equation} \label{eq:constants_Ass1}
    0< k_1,k_3 \ll k_2
\end{equation}
This peculiar structure of solutions has been observed early on as the Robertson model was widely used as a test problem for stiff numerical solvers, e.g., \cite[p.~3]{Hairer_1996}.
Up to our knowledge the Robertson model \eqref{equ:rob_reaction} has been investigated only numerically. No analytical results explaining the solution structure described above seem to be available. \comm{In this paper we provide a full asymptotic analysis of the Robertson model under the assumption $k_1, k_3 \ll k_2$, which covers the classical values \eqref{equ:rob_rates}. We rewrite the equations as a two-parameter singular perturbation problem in the rescaled small parameters $(\eps_1,\eps_2):=(k_1/k_2,k_3/k_2)$, which we then analyze using geometric singular perturbation theory (GSPT).}

\comm{Similar phenomena and solution structures occurr in many chemical reactions and more general classes of biological models. Due to the occurrence of variables and parameters of widely different orders of magnitude most of these models are multiscale in nature, see e.g \cite{Radulescu_Scaling_2021,segel1988validity,snowdensurvey} and references therein.
As a consequence of these multiple scales, some variables may vary little and can thus be treated as constants. Some parameters may have almost no effect and can be neglected. Some variables may rapidly approach a (quasi)equilibrium and can thus be slaved to other variables. Different mechanisms may dominate the dynamics in certain regions of phase- and/or parameter-space. Correspondingly, 
individual trajectories contain segments generated by
sequences of fast and slow processes on widely separated time scales. 
Identifying and analysing these regimes and the resulting decompositions into subsystems is crucial in the analysis of the dynamics. In suitably scaled variables the equations in the individual regimes often have the form of a slow-fast dynamical system and the dynamics is (locally) organized by a slow invariant manifold of lower dimension. The dynamics on the slow manifold is governed by a reduced model of lower dimension.\\
The most widely used realization of these ideas is the well known quasi steady state approximation (QSSA), which is used
to obtain lower-dimensional approximating models, i.e.
reactants involved in fast processes are eliminated by assuming that they are in equilibrium
\cite{segel1988validity,snowdensurvey}. }

A powerful \comm{mathematical} concept in explaining these phenomena are slow manifolds. The mathematical theory of slow manifolds and more general of slow-fast dynamical systems, known as geometric singular perturbation theory (GSPT), is well developed for ODEs depending singularly on one distinguished parameter $\varepsilon \ll 1$, see \cite{Jones_1995,Fenichel_1979,Krupa_2001_Extend, Kuehn_2015, wechselberger_geometric_2020} an the numerous references therein.
The origins of GSPT date back to the work of Fenichel \cite{Fenichel_1979}, where he introduced an invariant manifold  approach for singularly perturbed differential equations of the form 
\begin{equation} \label{eq:Fenichel_1}
    z'=H(z,\eps)
\end{equation}
with $z \in \R^k$, $k\geq 2$ and $\eps \ll 1$, see also \cite{wechselberger_geometric_2020} for a modern presentation. 
\comm{In this setting singular perturbation problems correspond to situations where the set
$$\mathcal{S}:=\{z \in \R^k: H(z,0) =0\}$$
contains manifolds of dimension $l \geq 1$. In many situations $\mathcal{S}$ is a manifold of dimension $l$, however, in other situations it is not a manifold in the strict sense due to the existence of singular points. Therefore $\mathcal{S}$ is denoted as the critical set. The critical set $\mathcal{S}$ corresponds to the equilibrium points of the $\eps=0$ limit of \eqref{eq:Fenichel_1} $z'=H(z,0)$. Under certain conditions manifolds $\Tilde{\mathcal{S}} \subset \mathcal{S}$ perturb to slow manifolds $\Tilde{\mathcal{S}}_{\eps}$ for $0<\eps \ll 1$ on which the dynamics is slow, while away from $\Tilde{\mathcal{S}}_{\eps}$ the dynamics is fast.}

An important special case of \eqref{eq:Fenichel_1} are slow-fast systems in standard form given by 
\begin{equation} \label{eq:Slow-fast_standard}
\begin{aligned}
    x'&=f(x,y,\eps)\\
    y'&=\eps g(x,y,\eps)
\end{aligned}
\end{equation}
with $x \in \R^m$, $y \in \R^n$ and $\eps \ll 1$, where differentiation is w.r.t the fast time $\tau$.
Systems of the form \eqref{eq:Slow-fast_standard} are called slow-fast in standard form, because as long as $f$ and $g$ are $\mathcal{O}(1)$ the dynamics of $x$ is fast compared to $y$, i.e. $x$ is the fast variable and $y$ the slow variable.

\begin{rem}
    It will turn out that for the analysis of the Robertson model both forms \eqref{eq:Fenichel_1} and \eqref{eq:Slow-fast_standard} are relevant. In the following explanation of the basic principles of GSPT, we will limit ourselves to the important special case \eqref{eq:Slow-fast_standard}.
\end{rem}

The $\eps=0$ limit problem of \eqref{eq:Slow-fast_standard}
\begin{equation} \label{eq:layer}
\begin{aligned}
    x'&=f(x,y,0)\\
    y'&=0
\end{aligned}
\end{equation}
is called \emph{layer problem}, which is used as an approximation of the fast dynamics. \comm{With a slight abuse of notation we denote the set of equilibria of \eqref{eq:layer}
$$\mathcal{S}=\{(x,y)^T \in \R^{m+n}:f(x,y,0)=0\},$$
as \emph{critical manifold}, despite the fact that $\mathcal{S}$ does not need to be a manifold in the strict sense (as pointed out above in the discussion of the general form \eqref{eq:Fenichel_1}).} 
By switching to the slow time $t=\eps \tau $ we may write system \eqref{eq:Slow-fast_standard} in the (for $\eps>0$) equivalent form 
\begin{equation} \label{eq:Slow-fast_standard_slow_time}
\begin{aligned}
    \eps \Dot{x}&=f(x,y,\eps)\\
    \Dot{y}&=g(x,y,\eps)
\end{aligned}
\end{equation}
where differentiation is w.r.t. the slow time $t$. The limit problem on the slow time scale 
\begin{equation} \label{eq:reduced}
\begin{aligned}
    0&=f(x,y,0)\\
    \Dot{y}&=g(x,y,0)
\end{aligned}
\end{equation}
is called \emph{reduced problem} and  is used as an approximation of the slow dynamics. Observe that the reduced problem is a dynamical system on the critical manifold $\mathcal{S}$. Parts of the critical manifold $\mathcal{S}$, where the Jacobian $\frac{\partial f}{\partial x}$ is regular, may be represented locally as graphs $x=h(y)$ by the implicit function theorem. The reduced flow on $\mathcal{S}$ is then given by 
$$\Dot{y}=g(h(y),y,0).$$

The goal of GSPT is to combine the dynamics of the two  simpler limiting systems \eqref{eq:layer}, and \eqref{eq:reduced} to understand the behaviour of \eqref{eq:Slow-fast_standard} for  $0 < \eps \ll 1$. 
In \cite{Fenichel_1979} Fenichel showed that if the Jacobian $\partial_x f$ is uniformly hyperbolic, the critical manifold $\mathcal{S}$ perturbs smoothly to a locally invariant slow manifold $\mathcal{S}_{\eps}$ which is $\mathcal{O}(\eps)$-close to $\mathcal{S}$, shares its stability properties with $\mathcal{S}$ and the slow flow on $\mathcal{S}_{\eps}$ converges to the reduced flow as $\eps \to 0$.

A major difficulty that remained in GSPT were non-hyperbolic points, i.e., points where at least one eigenvalue of the Jacobian $\partial_x f$ lies on the imaginary axis. Frequently these points are given by the singularities of the critical manifold. The problem remained open until the pioneering work of Dumortier and Roussarie \cite{dumortier1996canard} where they introduced the blow-up method, which was then
developed into a powerful tool in GSPT by Krupa and Szmolyan see \cite{Krupa_2001_Extend, krupa_extending_transcritical}. The main idea of the blow-up method is to first extend the state space by adding the trivial equation $\eps'=0$ and then introducing suitable weighted spherical coordinates to blow-up the singularity, e.g., a point to a sphere or a line to a cylinder. After dividing out a suitable power of the radial variable, less singular differential equations  are obtained
which often allow for a complete analysis with dynamical systems tools. By now the blow-up method has been widely used in the analysis of singularly perturbed differential equations, see e.g. \cite{Maesschalk_2008_Canard, Gucwa_2008, Jelbart_2021, Kosiuk_2016, Kosiuk_2011, Krupa_2001_Canard, Kuehn_Szmo_2015, Wechselberger_2001, Iuorio_MEMS_2019}.
It seems fair to say that  GSPT is very well developed for systems with a distinguished singular perturbation parameter $\varepsilon$ and that it has proven to be very useful in a large array of applications.

However, surprisingly little seems to be known in the case of systems
depending singularly on several small or large parameters,
e.g., chemical reactions with reaction rates $k_i$, $i=1,\ldots, p$ of widely differing orders of magnitude. 
\comm{Since this singular dependence on several parameters 
is rather the rule than the exception
for realistic chemical reactions and many other classes of biological models, it is important to develop methods suitable for 
the asymptotic analysis of such problems.\\
Our interest in the Robertson model comes from the fact, that it is a well known but fairly simple representative of this class of problems. 
We expect that the approach of this paper will be useful in
a wide variety of problems, e.g., various variants and extensions
of the basic Michaelis-Menten mechanism \cite{fraser1, roussel} and
models of the Belousov-Zhabotinskii reaction \cite{tyson1982}.
In ongoing work we are using this approach in the analysis of a five-variable model of the cell cycle \cite{deso,tyson1991} with singular dependence on three parameters.
}

An obvious and often used approach to apply \comm{asymptotic methods and in particular} GSPT to \comm{multi-parameter problems} 
is to reduce to the one-parameter case by identifying a suitable parameter $\varepsilon$ such that 
\begin{equation} \label{eq:epsilon_scaling}
    (k_1,\ldots,k_p)^T \sim (\varepsilon^{\alpha_1},\ldots,\varepsilon^{\alpha_p})^T, \, \alpha_i \in \Z,\, i=1,\ldots,p.
\end{equation}
A simple illustration of this approach (and its inherent arbitrariness) in the context of the Robertson model with the classical parameters \eqref{equ:rob_rates}
would be $\varepsilon=1/10$ which leads to $\alpha_1=2$, $\alpha_2=-7$, and $\alpha_3=-4$. This widely used approach, where parameters are restricted to a curve, can be very successful if good numerical values of the parameters are available, see, e.g., \cite{Radulescu_Scaling_2021, Jelbart_Calcium1_2022}. Unfortunately, this is often not the case \comm{and tools for qualitative analysis covering  wider ranges of parameters are needed }.

\comm{A powerful tool in finding significant scalings of parameters and variables in (polynomial) systems of differential equations 
based on Newton polyhedra and the associated power transformations
was developed by A. D. Bruno, see \cite{bruno2000power} and references therein. A related approach based on ideas from tropical geometry has been developed in \cite{Radulescu_Scaling_2021,deso}.
The connections between these approaches and our geometric approach will be explored in future work. Clearly, such problems can also be
treated by more conventional methods based on matched asymptotic
expansions.  An advantage of the geometric dynamical systems approach is that it provides (i) detailed insight into the underlying singular structures and (ii) leads to rigorous results
on the dynamics in appropriate ranges of the parameters.
}

Hence it is desirable to develop or adapt GSPT to problems depending singularly on several independent parameters $(\varepsilon_1,\ldots,\varepsilon_l)^T$, $l\geq 2$. Such problems are potentially more challenging since the singular behaviour and the multi-scale structure can vary significantly in a neighbourhood of the singular limit $(\varepsilon_1,\ldots,\varepsilon_l)^T=(0,\ldots,0)^T$. As a step towards a framework for multi-parameter singular perturbations of ODEs, we distinguish three different cases.
We expect that this classification is preliminary and not exhaustive, nevertheless
we feel it is useful as a first step.
For simplicity we phrase this classification for systems depending on two parameters, but it can be easily extended to systems depending on more parameters. 

{\bf Case 1:}
There exists an ordered sequence of time-scales, i.e., the system of differential equations has the form 
    \begin{equation} \label{eq:nested_timescales}
        \begin{aligned}
            \Dot{x}_1&=f_1(x,\eps_1,\eps_2)\\
            \Dot{x}_2&=\eps_1 f_2(x,\eps_1,\eps_2)\\
            \Dot{x}_3&=\eps_1 \eps_2 f_3(x,\eps_1,\eps_2),
        \end{aligned}
    \end{equation}
    with $0<\eps_1,\eps_2 \ll 1$, which is the three time-scale analogon to the slow-fast standard form \eqref{eq:Slow-fast_standard}. In this situation one can apply Fenichel theory iteratively to obtain a nested sequence of critical manifolds. This case is fairly well understood if the manifolds are normally hyperbolic, see \cite{Teixeira_nested_time_scale_2017}. If there are non-hyperbolic points, the situation can be more complicated, e.g., see the early influential paper \cite{ Krupa_Popovic_MMOS_2008} and the more recent \cite{Kuntz_Fold3D_2023}.
\par

In the two remaining cases, we consider more general systems in non-standard form, i.e., 
    \begin{equation} \label{eq:two_paramters_cases}
            \Dot{z}=H(z,\eps_1,\eps_2)
    \end{equation}
    with $0<\eps_1,\eps_2 \ll 1$.
    
{\bf Case 2:} 
The parameter $\eps_1$ is a classical singular perturbation parameter of \eqref{eq:two_paramters_cases} with corresponding critical manifold $\mathcal{S}(\eps_2)$ (depending on $\eps_2$) by standard Fenichel theory. The singular dependence of \eqref{eq:two_paramters_cases} on $\eps_2$ is caused by singularities of the critical manifold $\mathcal{S}(\eps_2)$ as $\eps_2 \to 0$, e.g., $\mathcal{S}(\eps_2)$ loses normal hyperbolicity, see \cite{Kosiuk_2011,Iuorio_MEMS_2019}.
\par

{\bf Case 3:}
Both parameters $\eps_1$ and $\eps_2$ act as singular perturbation parameters, leading to fundamentally different slow-fast structures in different regions of the parameter space.
\par

\comm{So far the analysis of problems corresponding to cases 2 and 3  has been carried out mostly
in the form of individual case studies}, e.g., see the very interesting work \cite{Dumortier_bogdanov_takens_2011} and also  \cite{carter_travelling_2023}. For more examples and an attempt to extract common features of existing results we refer to the recent review \cite{Kuehn_double_limits_2022} and the many references therein. 

\comm{The goal of this work is to make progress on adapting and extending GSPT to multi-parameter singular perturbation problems which is
the topic of the ongoing thesis project \cite{Baumgartner_thesis}. }
We give an asymptotic analysis of the Robertson model \eqref{equ:rob_reaction} under the assumption \eqref{eq:constants_Ass1}, which covers the classical choice \eqref{equ:rob_rates} in \cite{robertson_solution_1966}. It turns out that the Robertson model has features of case 2 and case 3, which shows that the above classification is not strict.
We view our analysis as a step in adapting GSPT to 
multi-parameter singular perturbation problems like \eqref{eq:two_paramters_cases} and also as a starting point for the analysis of similar problems depending on more than two parameters. First, we rewrite \eqref{equ:rob_reaction} as a two-parameter singular perturbation problem in the rescaled parameters 
$$(\varepsilon_1,\varepsilon_2)^T:=\Big(k_1 / k_2, k_3 / k_2 \Big)^T \in \R_+ \times \R_+, $$ 
varying in a neighbourhood of $(\eps_1,\eps_2)^T=(0,0)^T$.

Recall from the proof of Lemma \ref{lem:robertson_textbook}, that we can reduce the Robertson model to a planar dynamical system of the form \eqref{eq:rob_2D_intro}. By switching to the fast time scale $\tau=k_2t$ we obtain 
\begin{equation} \label{equ:rob_reaction_2D}
    \begin{aligned}
y'&=\varepsilon_1 (c-y-z) - y^2 - \varepsilon_2 yz\\
z'&= y^2,
\end{aligned}
\end{equation}
with initial value $(y_0,z_0)^T=(0,0)^T$, \enquote{ $ ' $ } denotes the derivative w.r.t. the fast time $\tau$, and $0<\eps_1$, $\eps_2 \ll 1$. System \eqref{equ:rob_reaction_2D} is now a planar multi-parameter singularly perturbed differential equation of the form \eqref{eq:two_paramters_cases}. Up to a reparametrization of time, system \eqref{equ:rob_reaction_2D} is equivalent to \eqref{equ:rob_reaction}, hence, we will perform our GSPT analysis based on  the planar system \eqref{equ:rob_reaction_2D}. 

\begin{rem}
    It follows from Lemma \ref{lem:robertson_textbook} that the solution of \eqref{equ:rob_reaction_2D} with initial value $(y_0,z_0)^T=(0,0)^T$
    converges to the equilibrium $Q=(0,c)^T$ for $\eps_1, \eps_2 > 0$. The linearization of \eqref{equ:rob_reaction_2D} at $Q$ has eigenvalues $\lambda_1=-\eps_1-\eps_2c$ and $\lambda_2=0$ with corresponding eigenvectors $v_1=(1,0)^T$ and $v_2=(\eps_1,-\eps_1-\eps_2c)^T$. Standard center manifold theory \cite{Guckenheimer_1983} implies that this solution converges to the equilibrium tangent to the \mbox{center-direction $v_2$.}
\end{rem}

\begin{figure}[t]
\centering
    \begin{tikzpicture}[scale=1]
			
        \draw[] (0,-0.5) -- (0,2);
        \draw[] (0,2) -- (0,3.5);
        \draw[] (-0.5,0) -- (2,0);
        \draw[] (2,0) -- (3.5,0);
        
        \draw[scale=1, domain=0:2.7, smooth, variable=\x, line width=0.5mm, dotted] plot ({\x}, {(1.3)*\x^(1/2)});
        \draw[scale=1, domain=0:1.4, smooth, variable=\x, line width=0.5mm, dotted] plot ({\x}, {(2.7)*\x^(1/2)});
        \draw[scale=1, domain=0:3.3, smooth, variable=\x, line width=0.5mm, dotted] plot ({\x}, {0.3*\x});
        \draw (83:2.5) node{$B_3$};
        \draw (60:2.5) node{$B_2$};
        \draw (30:2.5) node{$B_{12}$};
        \draw (7:2.5) node{$B_{11}$};
        \draw (-5:3) node{$\delta$};
        \draw (16:3.7) node{$C_1$};
        \draw (37:3.7) node{$C_2$};
        \draw (65:3.7) node{$C_3$};
        \draw (3.8,0) node{$\varepsilon_1$};
        \draw (0,3.8) node{$\varepsilon_2$};
        \filldraw[] (0,0) circle (2pt);
        \draw (3,0) arc (0:90:3);

			\end{tikzpicture}
   \caption{The four scaling regions $B_{11}$, $B_{12}$, $B_2$, and $B_3$.}
   \label{fig:scaling_in_parameter_space}
\end{figure}
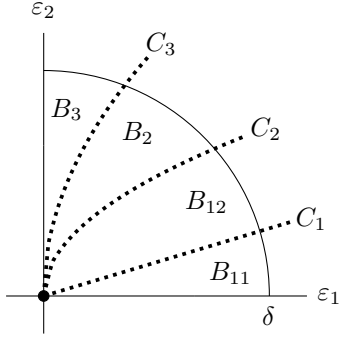

It turns out that for an asymptotic analysis, a small neighborhood of  $(\eps_1,\eps_2)^T=(0,0)^T$
must be divided into four regions corresponding to different singular limits and slow-fast structures in phase space, see Figure \ref{fig:scaling_in_parameter_space}. 
Our main result can be summarized as follows.

\begin{thm} \label{thm:main_summary}
    There exists $\delta>0$ such that the following holds in the $\delta$-neighbourhood 
    $$D_{\delta}:=\{(\varepsilon_1,\varepsilon_2)^T \in \R^2:\, \varepsilon_1 \geq 0,\, \varepsilon_2 \geq 0,\, \varepsilon_1^2+\varepsilon_2^2 \leq \delta \}$$ 
    of the origin in parameter space.
    \begin{enumerate}
        \item[1.] There exist constants $0<\beta_3<\beta_2$ and $\beta_1>0$ such that the curves $C_1=\{(\varepsilon_1,\varepsilon_2)^T \in \R^2:\varepsilon_1=\beta_1 \varepsilon_2\}$, $C_2=\{(\varepsilon_1,\varepsilon_2)^T \in \R^2:\varepsilon_1=\beta_2 \varepsilon_2^2\}$, and $C_3=\{(\varepsilon_1,\varepsilon_2)^T \in \R^2:\varepsilon_1=\beta_3 \varepsilon_2^2\}$ divide $D_{\delta}$ into four regions $B_{11}$, $B_{12}$, $B_2$, and $B_3$, see Figure \ref{fig:scaling_in_parameter_space}.
        \item[2.] In each of the regions  $B_{11}$, $B_{12}$, $B_2$, and $B_3$
        the problem \eqref{equ:rob_reaction_2D} has a different slow-fast structure each depending on a distinguished singular perturbation parameter. These structures become visible in suitable rescalings and blow-ups.
        \item[3.] For each of these regions  $B_{11}$, $B_{12}$, $B_2$, and $B_3$
        we identify a singular orbit $\gamma_0$ of a certain type connecting the initial value $O=(0,0)^T$ to the unique equilibrium $Q=(0,c)^T$ of \eqref{equ:rob_reaction_2D}. 
        \item[4.] In each of the regions  $B_{11}$, $B_{12}$, $B_2$, and $B_3$
        the orbit corresponding to the initial value approaches the corresponding singular orbit $\gamma_0$ in Hausdorff distance as $(\eps_1,\eps_2)^T \to (0,0)^T$ in the respective region, with error estimates depending on the sizes of $\varepsilon_1$, $\varepsilon_2$.
    \end{enumerate}
\end{thm}

\begin{rem}
(i)  By choosing slightly different constants $\beta_i$ the regions $B_{11}$, $B_{12}$, $B_2$, and $B_3$ can be viewed as overlapping. This implies that the multi-scale structure of the solution 
    changes in a smooth way for $\varepsilon_1$, $\varepsilon_2$ close to the curves $C_1$, $C_2$, and $C_3$. (ii)    Actually, Theorem \ref{thm:main_summary} holds for arbitrary constants  $0<\beta_3<\beta_2$ and $\beta_1>0$ if $\delta$ is chosen sufficiently small.
\end{rem}

\begin{figure}[h!]
    \centering
    \includegraphics[width=0.7\textwidth]{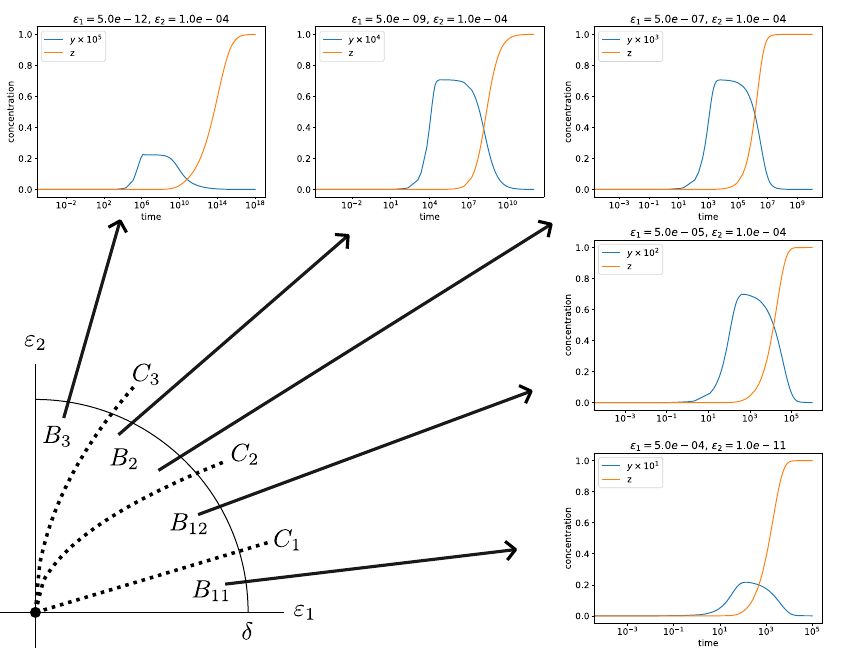}
    \caption{Numerical simulations of \eqref{equ:rob_reaction_2D} in different regions of the parameter space for $c=1$.}
    \label{fig:numerical_comparison}
\end{figure}

\comm{To give a first impression of the different slow-fast structures in these four scaling regions numerically computed solutions of \eqref{equ:rob_reaction_2D} are shown in Figure \ref{fig:numerical_comparison}. As one moves counterclockwise (increasing $\eps_2$ relative to $\eps_1$) one observes that the plateau value $y^{max}_{num}$ is shrinking and given by $y^{max}_{num}\approx 2.16\cdot 10^{-2}$, $y^{max}_{num}\approx 6.98\cdot 10^{-3}$, $y^{max}_{num}\approx 7.05\cdot 10^{-4}$, $y^{max}_{num}\approx 7.07\cdot 10^{-5}$, $y^{max}_{num}\approx 2.26\cdot 10^{-6}$ in $B_{11}$, $B_{12}$, $B_2$ (lower), $B_2$ (upper) and $B_3$, respectively. The structure of these solutions will be explained by relating them to the singular orbits $\gamma_0$ of Theorem \ref{thm:main_summary}. Note that the two solution profiles shown for region $B_2$ correspond to minor changes in the geometry of the underlying critical manifold, cf., the analysis in Section \ref{sec:Regime R2} summarized in Figure \ref{fig:K1_R2_rescaled_singular}.
In addition, we will show that the numerical values for $y^{max}_{num}$ fit well with corresponding values obtained by our analysis, see Section \ref{sec:Summary}.  
}

Our analysis and proofs are based on suitable blow-ups of the origin in parameter space which combined with blow-ups in phase space reveal the underlying slow-fast structures in the regions $B_{11}$, $B_{12}$, $B_2$, and $B_3$. We are confident that this approach can also be useful in the analysis of systems with more than two singular perturbation parameters.

The rest of the paper is organized as follows: In a first step it is convenient to
blow up the origin in parameter space $(\eps_1,\eps_2)^T =(0,0)^T$ in a suitable way. This is done in 
Section~\ref{sec:Parameter blow-up}. Loosely speaking this allows 
to apply GSPT with the radial parameter as a distinguished singular perturbation parameter.
In Section \ref{sec:Regime R2} we carry out the rather straightforward
GSPT analysis for region $B_2$. 
The slow fast-structures corresponding to the regions $B_{11}$, $B_{12}$ and $B_3$ are more complicated and require additional blow-ups. The analysis of these cases is carried out in
in Section \ref{sec:Regime R1} and \ref{sec:Regime R3}, respectively. We end with a conclusion and outlook.

\section{Structure of Parameter Space} \label{sec:Parameter blow-up}

The goal in singularly perturbed systems with a single parameter $\eps \ll 1$ is to prove statements which hold for $\eps \in (0,\Hat{\eps}]$ for some $\Hat{\eps}>0$. In system \eqref{equ:rob_reaction_2D} we are now dealing with a two-parameter problem in $\eps_1, \, \eps_2 \ll 1$, hence we need to prove results which hold in a small neighbourhood of the origin in the parameter space $\R^2_+$.

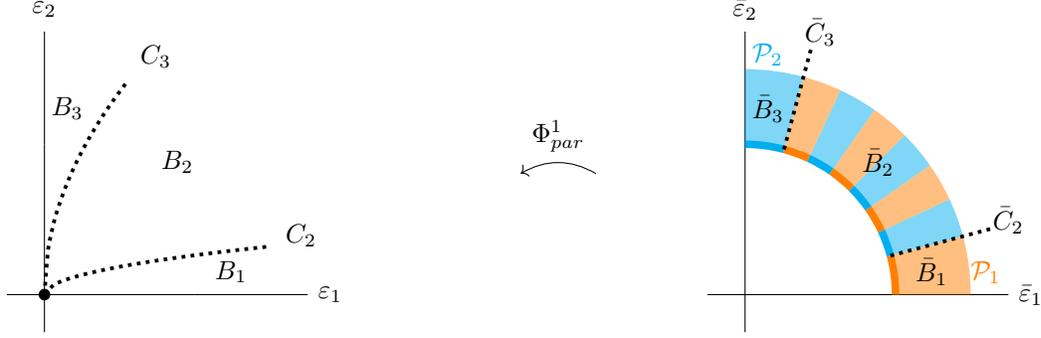
\begin{figure}
\centering
\begin{subfigure}{0.4\textwidth}
    \begin{tikzpicture}[scale=1]
			
        \draw[] (0,-0.5) -- (0,2);
        \draw[] (0,2) -- (0,3.5);
        \draw[] (-0.5,0) -- (2,0);
        \draw[] (2,0) -- (3.5,0);
        
        \draw[scale=1, domain=0:3, smooth, variable=\x, line width=0.5mm, dotted] plot ({\x}, {(1/2.7)*\x^(1/2)});
        \draw[scale=1, domain=0:1.1, smooth, variable=\x, line width=0.5mm, dotted] plot ({\x}, {(2.7)*\x^(1/2)});
        \draw (83:2.5) node{$B_3$};
        \draw (45:2.5) node{$B_2$};
        \draw (7:2.5) node{$B_1$};
        \draw (13:3.5) node{$C_2$};
        \draw (65:3.5) node{$C_3$};
        \draw (3.8,0) node{$\varepsilon_1$};
        \draw (0,3.8) node{$\varepsilon_2$};
        \filldraw[] (0,0) circle (2pt);

			\end{tikzpicture}
\end{subfigure}
\begin{subfigure}{0.15\textwidth}
\begin{tikzpicture}
       \path[->]
    (1,2) edge[bend right] node [left] {} (0,2);
    \draw (0.5,2.5) node{$\Phi_{par}^1$};
    \draw (0,0) node{};
\end{tikzpicture}
\end{subfigure}
\begin{subfigure}{0.4\textwidth}
    \begin{tikzpicture}[scale=1]
        \draw[] (0,-0.5) -- (0,3.5);
        \draw[] (-0.5,0) -- (3.5,0);
        
        \draw[orange,line width=1mm] (2, 0) arc (0:15:2);
        \draw[cyan,line width=1mm] (0,2) arc (90:75:2);
        
        \draw[orange,line width=1mm] (75:2) arc (75:65:2);
        \draw[cyan,line width=1mm] (65:2) arc (65:55:2);
        \draw[orange,line width=1mm] (55:2) arc (55:45:2);
        \draw[cyan,line width=1mm] (45:2) arc (45:35:2);
        \draw[orange,line width=1mm] (35:2) arc (35:25:2);
        \draw[cyan,line width=1mm] (25:2) arc (25:15:2);
        \fill[orange,opacity=.5]
            (0:2) -- (0:3) -- (0:3) arc(0:15:3)-- (15:2) --cycle;
        \fill[cyan,opacity=.5]
            (15:2) -- (15:3) -- (15:3) arc(15:25:3)-- (25:2) --cycle;
        \fill[orange,opacity=.5]
            (25:2) -- (25:3) -- (25:3) arc(25:35:3)-- (35:2) --cycle;
        \fill[cyan,opacity=.5]
            (35:2) -- (35:3) -- (35:3) arc(35:45:3)-- (45:2) --cycle;
            \fill[orange,opacity=.5]
            (45:2) -- (45:3) -- (45:3) arc(45:55:3)-- (55:2) --cycle;
            \fill[cyan,opacity=.5]
            (55:2) -- (55:3) -- (55:3) arc(55:65:3)-- (65:2) --cycle;
            \fill[orange,opacity=.5]
            (65:2) -- (65:3) -- (65:3) arc(65:75:3)-- (75:2) --cycle;
            \fill[cyan,opacity=.5]
            (75:2) -- (75:3) -- (75:3) arc(75:90:3)-- (90:2) --cycle;
            \draw[line width=0.5mm,dotted] (15:2) -- (15:3.4);
            \draw[line width=0.5mm,dotted] (75:2) -- (75:3.4);
        \draw (83:2.5) node{$\Bar{B}_3$};
        \draw (45:2.5) node{$\Bar{B}_2$};
        \draw (7:2.5) node{$\Bar{B}_1$};
        \draw (3.8,0) node{$\Bar{\varepsilon}_1$};
        \draw (0,3.8) node{$\Bar{\varepsilon}_2$};
        \draw[orange,line width=1mm] (3.2,0.3) node{$\mathcal{P}_1$};
        \draw[cyan,line width=1mm] (0.3,3.2) node{$\mathcal{P}_2$};
        \draw[] (3.5,1) node{$\Bar{C}_2$};
        \draw[] (1,3.5) node{$\Bar{C}_3$};
    
        

			\end{tikzpicture}
\end{subfigure}
			\caption{Parameter blow-up $\Phi_{par}^1$ of the origin and charts $\mathcal{P}_1$ (orange) and $\mathcal{P}_2$ (blue) }
			\label{fig:parameter_blow-up}
\end{figure}

As a first step, it is instructive to look at the three limiting problems of \eqref{equ:rob_reaction_2D}:
\begin{enumerate}
    \item[1)] $\varepsilon_2=0,\,\varepsilon_1>0:$ There exists a unique equilibrium given by $(y,z)^T=(0,c)^T$.
The linearization at the equilibrium has one negative and one vanishing eigenvalue, thus center manifold theory can be applied there.
        \item[2)] $\varepsilon_1=0,\,\varepsilon_2>0:$ The line $y=0$ consists of equilibria. 
    The line of equilibria $\{(0,z)^T, \, z >0\}$ is  attracting for $\eps_2 >0$. 
    The origin $(y,z)^T=(0,0)^T$ is more degenerate, i.e. the corresponding linearization has a double zero
    eigenvalue.
    \item[3)] $\varepsilon_1=0=\varepsilon_2:$ The line $y=0$ consists of degenerate equilibria, 
i.e. the corresponding linearizations have a double zero
    eigenvalue.
    
\end{enumerate}
The three cases above are qualitatively quite different, ranging from a unique equilibrium, which can be analysed by center manifold reduction, to a very degenerate line of nilpotent equilibria. This indicates that in the double limit we should expect that the relative sizes of $\varepsilon_1$ and $\varepsilon_2$ have a significant influence on the detailed dynamics and asymptotics. 
It turns out that this is indeed the case and parameter space must be divided into three regions $B_1$, $B_2$ and $B_3$ where 
$$\varepsilon_2^2 \ll \varepsilon_1, \quad \varepsilon_1 \approx \varepsilon_2^2, \quad \varepsilon_1 \ll \varepsilon_2^2,$$
respectively. To be precise we define the curves 
\begin{equation} \label{eq:parameter_blow-up_curves}
    C_2:=\{(\varepsilon_1,\varepsilon_2)^T \in \R^2:\varepsilon_1=\beta_2 \varepsilon_2^2\},\quad  C_3=:\{(\varepsilon_1,\varepsilon_2)^T \in \R^2:\varepsilon_1=\beta_3 \varepsilon_2^2\}
\end{equation}
for $0<\beta_3<\beta_2$ and the regions 
\begin{align} 
   B_1&= \{(\varepsilon_1,\varepsilon_2)^T \in \R^2:\varepsilon_1 > \beta_2\varepsilon_2^2 \} \label{equ:parameter_regions_R1_before_blowup}\\ 
   B_2&=\{(\varepsilon_1,\varepsilon_2)^T \in \R^2:\beta_3\varepsilon_2^2 \leq \varepsilon_1 \leq \beta_2\varepsilon_2^2 \} \label{equ:parameter_regions_R2_before_blowup}\\
   B_3&=\{(\varepsilon_1,\varepsilon_2)^T \in \R^2:\varepsilon_1 < \beta_3\varepsilon_2^2 \},\label{equ:parameter_regions_R3_before_blowup}
\end{align}
see Figure \ref{fig:parameter_blow-up} (left).

To separate the curves $C_2$ and $C_3$ in a neighbourhood of the origin we perform a non-homogeneous blow-up transformation. It turns out that this allows for a 
GSPT analysis in Region $B_2$, by using the radial parameter as singular perturbation parameter.

The blow-up map respecting the scaling properties of the curves $C_2$ and $C_3$ is
\begin{equation} \label{equ:parameter_blow-up}
 \begin{aligned}
         \Phi_{par}^1: [0,\infty) \times \mathbb{S}^1 &\to \R^2\\
         (r,\Bar{\varepsilon}_1,\Bar{\varepsilon}_2) &\mapsto \begin{cases}
             \eps_1=r^2\Bar{\varepsilon}_1\\
             \eps_2=r\Bar{\varepsilon}_2,
         \end{cases}
 \end{aligned}
\end{equation}
where we naturally restrict ourselves to the meaningful parameter space $\Bar{\varepsilon}_1,\Bar{\varepsilon}_2\geq 0$.   The preimage of the origin under $\Phi_{par}^1$ is the quarter circle ($r=0$), which implies that $\Phi_{par}^1$ is not injective for $r=0$. Away from the origin the blow-up map $\Phi_{par}^1$ is a diffeomorphism. 
In the blown-up parameter space the quadratic curves $C_2$ and $C_3$ correspond to well separated straight lines $\Bar{C}_2$ and $\Bar{C}_3$ given by
\begin{align} 
   \Bar{C}_2&=\Big\{(r,\Bar{\varepsilon}_1,\Bar{\varepsilon}_2)^T \in [0,\infty) \times \mathbb{S}^1:\Bar{\varepsilon}_1 = -1/2 \beta_2+\sqrt{1/4\beta_2^2+1}\Big\} \label{equ:parameter_regions_R1}\\ 
   \Bar{C}_3&=\Big\{(r,\Bar{\varepsilon}_1,\Bar{\varepsilon}_2)^T \in [0,\infty) \times \mathbb{S}^1: \Bar{\varepsilon}_1=-1/2 \beta_3+\sqrt{1/4\beta_3^2+1} \Big\} \label{equ:parameter_regions_R2}
\end{align}
with $\beta_3<\beta_2$ from \eqref{eq:parameter_blow-up_curves}, respectively, see Figure \ref{fig:parameter_blow-up}.
The regions $B_1$, $B_2$, and $B_3$ correspond to $\Bar{B}_1$, $\Bar{B}_2$, and $\Bar{B}_3$ in the obvious way. The size of the constants $\beta_2$ and $\beta_3$ determines the size of the regions $\Bar{B}_1$, $\Bar{B}_2$, and $\Bar{B}_3$. 
For $\beta_2 \to \infty$ the line $\Bar{C}_2$ approaches the $\Bar{\eps}_1$-axis, 
similarly the line $\Bar{C}_3$   approaches the $\Bar{\eps}_2$-axis as  $\beta_3 \to 0$.
\begin{rem}
    The choice of the constants $\beta_2$ and $\beta_3$ determines the size of the neighbourhood in which our GSPT analysis is valid. However, for arbitrary constants $0<\beta_3<\beta_2$ we can always find such a sufficiently small neighbourhood.
\end{rem}

It is natural to perform the remaining analysis in directional charts $\mathcal{P}_1$ and $\mathcal{P}_2$ corresponding to the directions $\Bar{\varepsilon}_1 = 1$ and $\Bar{\varepsilon}_2 =1$,
respectively. In these charts the blow-up transformation has the form 
\begin{align} 
    \mathcal{P}_1&: \varepsilon_1=r^2,\, \varepsilon_2=r\tilde{\varepsilon}_2 \label{equ:parameter_blow-up-chart_K1}\\
    \mathcal{P}_2&: \varepsilon_1=r^2\tilde{\varepsilon}_1,\, \varepsilon_2=r,\label{equ:parameter_blow-up-chart_K3}
\end{align}
respectively. Chart $\mathcal{P}_1$ covers the regions  $\Bar{B}_1$ and $\Bar{B}_2$,
while $\mathcal{P}_2$ covers the regions  $\Bar{B}_2$ and $\Bar{B}_3$, 
see Figure \ref{fig:parameter_blow-up} where the regions covered by charts 
$\mathcal{P}_1$ and $\mathcal{P}_2$ are shown in orange and blue, respectively.
The alternating colors in region $\Bar{B}_2$ indicate that this region is covered by both charts. 

The regions $\Bar{B}_1$ and $\Bar{B}_2$ in chart $\mathcal{P}_1$ are given by $0\leq \tilde{\varepsilon}_2 < \sqrt{\frac{1}{\beta_2}}$ and $\sqrt{\frac{1}{\beta_2}} \leq \tilde{\varepsilon}_2 \leq \sqrt{\frac{1}{\beta_3}}$, respectively. For the analysis in region $\Bar{B}_3$ its description in chart $\mathcal{P}_2$, i.e., $0 \leq \tilde{\varepsilon}_1 < \beta_3$, will be relevant.

We start with the analysis in region $\Bar{B}_2$, which is the simplest case and covers the 
slow-fast structure corresponding to the classical parameters \eqref{equ:rob_rates}.
The regions $\Bar{B}_1$ and $\Bar{B}_3$ correspond to more degenerate cases and somewhat more complicated slow-fast structures, which we will treat afterwards.

\section{Analysis in Region $B_2$} \label{sec:Regime R2}
The analysis in region $\Bar{B}_2$ can be carried out 
in any of the two charts $\mathcal{P}_i$, $i=1,2$, we choose to work in chart $\mathcal{P}_1$. Inserting \eqref{equ:parameter_blow-up-chart_K1} into \eqref{equ:rob_reaction_2D}, we obtain a slow-fast system in non-standard form 
\begin{equation} \label{equ:rob_reaction_2D_K1}
    \begin{aligned}
y'&=r^2 (c-y-z) - y^2 - r \tilde{\eps}_2 yz\\
z'&= y^2,
\end{aligned}
\end{equation}
where $r$, $\tilde{\eps}_2 \in \R_{\geq 0}$. It will be important that in region $\Bar{B}_2$ 
we have $\tilde{\varepsilon}_2 \in \Big[\sqrt{\frac{1}{\beta_2}}, \sqrt{\frac{1}{\beta_3}}\Big]$.
This avoids degeneracies occurring as  $\tilde{\eps}_2 \to 0$  or  $\tilde{\eps}_2 \to \infty$ 
which are treated in the analysis of regions $\Bar{B}_1$ and $\Bar{B}_3$. 
For better readability we are dropping the \enquote{ $\Tilde{ }$ } in the following.

 \begin{figure}
 \centering
			\begin{tikzpicture}[scale=1.5]
        \draw[->-=.7,black!50!green] (0.5,-1) -- (0,-0.5);
        \draw[->-=.8,black!50!green] (0,-0.5) -- (-1,0.5);
        \draw[->-=.7,black!50!green] (1.5,-0.75) -- (0,0.75);
        \draw[->-=.6,black!50!green] (0,0.75) -- (-1,1.75);
        \draw[->-=.6,black!50!green] (1.5,0.5) -- (0,2);
        \draw[->-=.6,black!50!green] (0,2) -- (-0.5,2.5);
        
        \draw[blue] (0,-1) -- (0,2.5);
        \draw[] (-1,0) -- (2.5,0);
        \draw (2.7,0) node{$y$};
        \draw (0,2.7) node{$z$};
        \draw[] (-0.15,0.15) node[] {O};
        \draw[] (-0.15,1.65) node[] {Q};
        \filldraw[teal] (0,0) circle (2pt);
        \filldraw[black] (0,1.5) circle (2pt);
			\end{tikzpicture}
\caption{Phase portrait of the layer problem \eqref{eq:degenerate_limit}.}
\label{fig:rob_reaction_2d_k1=k3=0}
\end{figure}
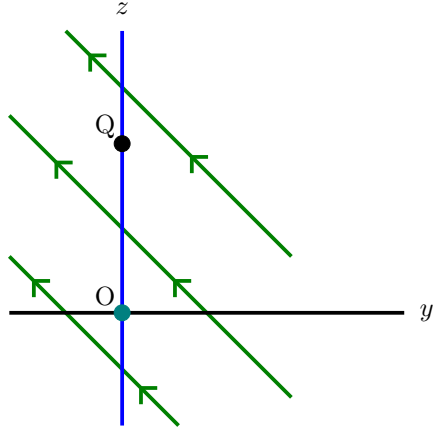

In system \eqref{equ:rob_reaction_2D_K1}, the parameter $r$ is the slow-fast parameter. The layer problem ($r=0$) has the simple form 
\begin{equation}
    \begin{aligned} \label{eq:degenerate_limit}
    y'&=-y^2\\
    z'&=y^2,
\end{aligned}
\end{equation}

 which obviously coincides with the limit problem $\eps_1=\eps_2=0$ of \eqref{equ:rob_reaction_2D}. System \eqref{eq:degenerate_limit} is explicitly solvable and its orbits are straight lines with slope $-1$, i.e.,
 $$z=-y+s, \quad  s \in \R.$$ 
As mentioned before, $y=0$ is a line of nilpotent equilibria which attracts all orbits with $y(0)>0$ in forward time and attracts all orbits with $y(0)<0$ in backward time.  As a consequence of this degeneracy, solutions are very sensitive to perturbations around $y=0$. Note that the initial value $O=(0,0)^T$ and also the unique equilibrium $Q=(0,c)^T$ of \eqref{equ:rob_reaction_2D_K1} lie on the line of equilibria $y=0$ represented by a teal and black dot in Figure \ref{fig:rob_reaction_2d_k1=k3=0}, respectively.
 
 In terms of slow-fast systems the critical manifold $\mathcal{S}$ is given by 
 $$\mathcal{S}=\{(y,z)^T \in \R^2:y=0\},$$ 
 which is not normally hyperbolic (which is indicated by green simple arrows in Figure \ref{fig:rob_reaction_2d_k1=k3=0}). Due to the lack of normal hyperbolicity Fenichel theory is not applicable. We resolve this degeneracy by rescaling the variable $y$ with
\begin{equation} \label{eq:rescaling_R2}
    y=r\Tilde{y}.
\end{equation}
Inserting \eqref{eq:rescaling_R2} into \eqref{equ:rob_reaction_2D_K1} gives 
\begin{equation} \label{eq:rob_reaction_2D_R2_singular}
    \begin{aligned}
\Tilde{y}'&=r(c-r \Tilde{y}-z - \Tilde{y}^2 -  \eps_2 \Tilde{y}z)\\
z'&= r^2 \Tilde{y}^2.
\end{aligned}
\end{equation}
For $r=0$ this vector field vanishes identically, thus we desingularize the system by dividing out a factor $r$, which can be viewed as transforming to  a slower timescale. Clearly this does not change the orbits of the system. \comm{This leads to
\begin{equation} \label{eq:rob_reaction_2D_R2_rescaled_tilde}
    \begin{aligned}
\Tilde{y}'&=c-r \Tilde{y}-z - \Tilde{y}^2 -  \eps_2 \Tilde{y}z\\
z'&= r \Tilde{y}^2.
\end{aligned}
\end{equation}
System \eqref{eq:rob_reaction_2D_R2_rescaled_tilde} is of standard slow-fast form w.r.t. the singular perturbation parameter $r$. To simplify the notation in the following computations we drop the \enquote{ $\Tilde{ }$ } and obtain 
\begin{equation} \label{eq:rob_reaction_2D_R2_rescaled}
    \begin{aligned}
y'&=c-r y-z - y^2 -  \eps_2 yz\\
z'&= r y^2.
\end{aligned}
\end{equation}
which will be the starting point for the analysis throughout the rest of the analysis in Section \ref{sec:Regime R2} and Section \ref{sec:Regime R1}.}

\begin{figure}[h!]
    \begin{subfigure}{0.3\textwidth}
        \centering
			\begin{tikzpicture}[
			scale=1.5]
			
        \draw[->-=.45,scale=1, domain=0:0.8,samples=500, variable=\x, red] plot ({\x}, {0.1-(1/10)*(0.8-\x)^(1/8)});
        \draw[scale=1,red] (0.8,0.06) -- (0.8,0.11);
        \draw[scale=1, domain=0.8:0.67, variable=\x, red] plot ({\x}, {0.1+(1/1.2)*(0.8-\x)^(1/2)});
        \draw[->-=.45,scale=1, domain=0.67:0, variable=\x, red] plot ({\x}, {(1-\x*\x)/(1+0.55*\x)});
        \draw[] (0.6,0) node[above] {\textcolor{red}{$\gamma_r$}};
\draw[] (0,-1) -- (0,1.5,0) node[above]{\textcolor{black}{$z$}};
        \draw[->>-=.7,black!50!green] (-1,0) -- (-1.5,0);
        \draw[->>-=.8,black!50!green] (-1,0) -- (1,0);
        \draw[->>-=.7,black!50!green] (1.5,0) node[right]{\textcolor{black}{$y$}} -- (1,0) ;
        \draw[->>-=.5,black!50!green] (1.2,0.5) -- (0.8,0.5);
        \draw[->>-=.5,black!50!green] (0,0.5) -- (0.8,0.5);
        \draw[->>-=.95,black!50!green] (0,0.5) -- (-1.2,0.5);
        \draw[->>-=.2,black!50!green] (1,1) -- (0,1);
        \draw[->>-=.5,black!50!green] (-0.4,1) -- (0,1);
        \draw[->>-=.9,black!50!green] (-0.4,1) -- (-1,1);
        \draw[->-=.45,scale=1, domain=1.5:-0.27, variable=\x, blue] plot ({\x}, {(1-\x*\x)/(1+0.5*\x)});
        \draw[scale=1, domain=-1.15:-0.27, variable=\x, blue] plot ({\x}, {(1-\x*\x)/(1+0.5*\x)});
        \filldraw[teal] (0,0) circle (2pt);
        \filldraw[black] (0,1) circle (2pt);
        \draw[] (-0.15,-0.15) node[] {O};
        \draw[] (0.15,1.15) node[] {Q};
        \draw[] (0.3,0.05) node[below] {\textcolor{black!50!green}{$\gamma_0^f$}};
        \draw[] (0.6,0.5) node[above] {\textcolor{blue}{$\gamma_0^s$}};
        \draw[-stealth,very thick,black] (0.5,-0.5) node[right,black] {\small $y^{max}$} -- ++ (0.45,0.45) ;
			\end{tikzpicture}
			\caption{$\frac{1}{\eps_2^2}>c$.}
			\label{fig:K1_R2_rescaled_singular1}
    \end{subfigure}
    \hfill
    \begin{subfigure}{0.3\textwidth}
        \centering
			\begin{tikzpicture}[scale=1.5]
			
        \draw[->-=.45,scale=1, domain=0:0.78,samples=500, variable=\x, red] plot ({\x}, {0.1-(1/10)*(0.8-\x)^(1/8)});
        \draw[->-=.5,scale=1, domain=0.8:0, variable=\x, red] plot ({\x}, {1-(0.86/0.8)*\x});
        \draw[smooth,red] plot coordinates { (0.78,0.04) (0.79,0.06) (0.81,0.1) (0.81,0.12) (0.8,0.14)  };
        \draw[] (0.6,0) node[above] {\textcolor{red}{$\gamma_r$}};
\draw[] (0,-1) -- (0,1.5,0) node[above]{\textcolor{black}{$z$}};
        \draw[->>-=.8,black!50!green] (-0.5,0) -- (1,0);
        \draw[->>-=.7,black!50!green] (1.5,0) node[right]{\textcolor{black}{$y$}} -- (1,0) ;
        \draw[->>-=.5,black!50!green] (1.2,0.5) -- (0.8,0.5);
        \draw[->>-=.3,black!50!green] (-0.5,0.5) -- (0.8,0.5);
        \draw[->>-=.9,black!50!green] (1,1) -- (0.8,1);
        \draw[->>-=.2,black!50!green] (-0.5,1) -- (0.8,1);
        \draw[->-=.4,scale=1, domain=1.5:-0.27, variable=\x, blue] plot ({\x}, {1-\x});
        \filldraw[teal] (0,0) circle (2pt);
        \filldraw[black] (0,1) circle (2pt);
        \draw[] (0.3,0.05) node[below] {\textcolor{black!50!green}{$\gamma_0^f$}};
        \draw[] (0.5,0.5) node[above] {\textcolor{blue}{$\gamma_0^s$}};
        \draw[] (-0.15,-0.15) node[] {O};
        \draw[] (0.15,1.15) node[] {Q};
        \draw[-stealth,very thick,black] (0.5,-0.5) node[right,black] {\small $y^{max}$} -- ++ (0.45,0.45) ;
        
        ;
			
			\end{tikzpicture}
			\caption{$\frac{1}{\eps_2^2}=c$.}
			\label{fig:K1_R2_rescaled_singular3}
    \end{subfigure}
    \hfill
    \begin{subfigure}{0.3\textwidth}
        \centering
			\begin{tikzpicture}[scale=1.5]
			
        \draw[->-=.45,scale=1, domain=0:0.78,samples=500, variable=\x, red] plot ({\x}, {0.1-(1/10)*(0.8-\x)^(1/8)});
        \draw[->-=.5,scale=1, domain=0.8:0.2, variable=\x, red] plot ({\x}, {(1-(\x+0.03)*(\x+0.03))/(1+1.5*(\x+0.03))});
        \draw[smooth,red] plot coordinates { (0.78,0.04) (0.79,0.06) (0.81,0.1) (0.81,0.12) (0.8,0.14) };
        \draw[scale=1, domain=0.2:0, variable=\x, red] plot ({\x}, {1/(1+2.1*\x)});
        \draw[] (0.6,0) node[above] {\textcolor{red}{$\gamma_r$}};
\draw[] (0,-1) -- (0,1.5,0) node[above]{\textcolor{black}{$z$}};
        \draw[->>-=.8,black!50!green] (-0.5,0) -- (1,0);
        \draw[->>-=.7,black!50!green] (1.5,0) node[right]{\textcolor{black}{$y$}} -- (1,0) ;
        \draw[->>-=.5,black!50!green] (1.2,0.5) -- (0.8,0.5);
        \draw[->>-=.3,black!50!green] (-0.5,0.5) -- (0.8,0.5);
        \draw[->>-=.9,black!50!green] (1,1) -- (0.8,1);
        \draw[->>-=.2,black!50!green] (-0.5,1) -- (0.8,1);
        \draw[->-=.4,scale=1, domain=1.5:-0.27, variable=\x, blue] plot ({\x}, {(1-\x*\x)/(1+1.5*\x)});
        \filldraw[teal] (0,0) circle (2pt);
        \filldraw[black] (0,1) circle (2pt);
        \draw[] (0.3,0.05) node[below] {\textcolor{black!50!green}{$\gamma_0^f$}};
        \draw[] (0.45,0.5) node[above] {\textcolor{blue}{$\gamma_0^s$}};
        \draw[] (-0.15,-0.15) node[] {O};
        \draw[] (0.15,1.15) node[] {Q};
        \draw[-stealth,very thick,black] (0.5,-0.5) node[right,black] {\small $y^{max}$} -- ++ (0.45,0.45) ;
        
        ;
			
			\end{tikzpicture}
			\caption{$\frac{1}{\eps_2^2}<c$.}
			\label{fig:K1_R2_rescaled_singular2}
    \end{subfigure}
    \caption{Singular orbit structure (green and blue) of \eqref{eq:rob_reaction_2D_R2_rescaled} and genuine orbit (red) connecting $O$ and $Q$ for $0<r \ll 1$.}
    \label{fig:K1_R2_rescaled_singular}
\end{figure}
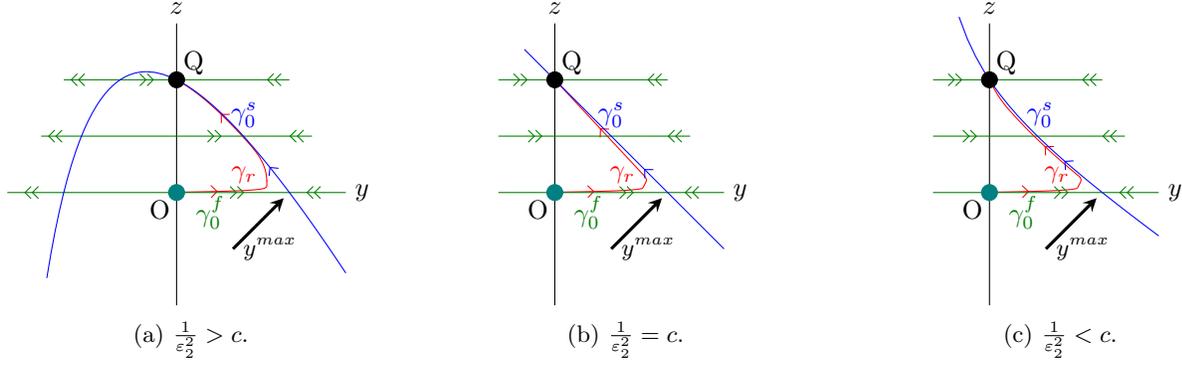

For $r=0$ we obtain the layer problem
\begin{equation} \label{eq:rob_reaction_2D_R2_rescaled_singular}
    \begin{aligned}
y'&=c-z - y^2 -  \eps_2 yz\\
z'&=0.
\end{aligned}
\end{equation}
The critical manifold is $\mathcal{S}=\{(x,y)^T \in \R^2: c-z-y^2-\eps_2yz=0 \}$. 
In the following we focus on the part of $\mathcal{S}$ in the half plane $y\geq 0$, denoted by $\mathcal{S}^a$, which is normally attracting for $\varepsilon_2 \in \Big[\sqrt{\frac{1}{\beta_2}}, \sqrt{\frac{1}{\beta_3}}\Big]$ and
can be described as a graph 
\begin{equation}
    z=\frac{c-y^2}{1+\eps_2y}.
\end{equation} 
The hyperbolicity of $\mathcal{S}^a$ follows since the eigenvalue of the corresponding linearization of \eqref{eq:rob_reaction_2D_R2_rescaled_singular} is $\lambda_1=-2y-\eps_2z < 0$. 
Note that $\mathcal{S}^a$ intersects the positive $y$-axis at $y = y^{max} :=\sqrt{c}$, see Figure \ref{fig:K1_R2_rescaled_singular}.

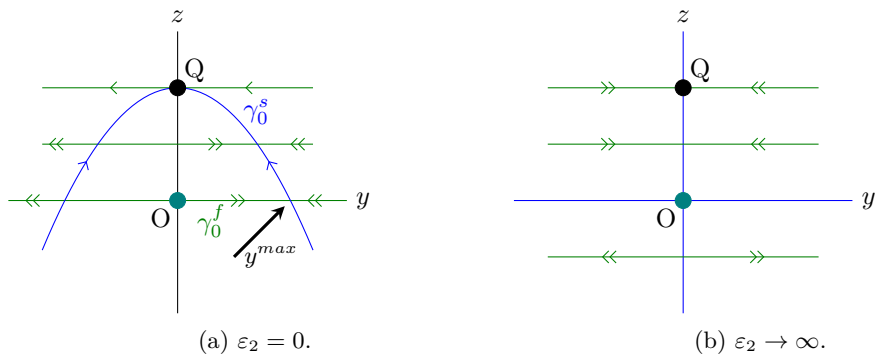
\begin{figure}
\centering
\begin{subfigure}{0.4\textwidth}
    \begin{tikzpicture}[
			scale=1.5]
			
\draw[] (0,-1) -- (0,1.5,0) node[above]{\textcolor{black}{$z$}};
        \draw[->-=.5,black!50!green] (1.2,1) -- (0,1);
        \draw[->-=.5,black!50!green] (0,1) -- (-1.2,1);
        \draw[->>-=.7,black!50!green] (-1,0) -- (-1.5,0);
        \draw[->>-=.8,black!50!green] (-1,0) -- (1,0);
        \draw[->>-=.7,black!50!green] (1.5,0) node[right]{\textcolor{black}{$y$}} -- (1,0) ;
        \draw[->>-=.5,black!50!green] (1.2,0.5) -- (0.8,0.5);
        \draw[->>-=.5,black!50!green] (0,0.5) -- (0.8,0.5);
        \draw[->>-=.95,black!50!green] (0,0.5) -- (-1.2,0.5);
        \draw[->-=.45,scale=1, domain=1.2:0, variable=\x, blue] plot ({\x}, {1-\x*\x});
        \draw[->-=.45,scale=1, domain=-1.2:0, variable=\x, blue] plot ({\x}, {1-\x*\x});
        \filldraw[teal] (0,0) circle (2pt);
        \filldraw[black] (0,1) circle (2pt);
        \draw[] (-0.15,-0.15) node[] {O};
        \draw[] (0.15,1.15) node[] {Q};
        \draw[] (0.3,0.05) node[below] {\textcolor{black!50!green}{$\gamma_0^f$}};
        \draw[] (0.7,0.6) node[above] {\textcolor{blue}{$\gamma_0^s$}};
        \draw[-stealth,very thick,black] (0.5,-0.5) node[right,black] {\small $y^{max}$} -- ++ (0.45,0.45) ;
			\end{tikzpicture}
   \caption{$\eps_2=0$.}
			\label{fig:K1_R1_rescaled_singular_eps=0}
\end{subfigure}
\begin{subfigure}{0.4\textwidth}
    \begin{tikzpicture}[
			scale=1.5]
			
\draw[blue] (0,-1) -- (0,1.5,0) node[above]{\textcolor{black}{$z$}};
        \draw[blue] (1.5,0) node[right]{\textcolor{black}{$y$}} -- (-1.5,0) ;
        \draw[->>-=.5,black!50!green] (1.2,0.5) -- (0,0.5);
        \draw[->>-=.5,black!50!green] (-1.2,0.5) -- (0,0.5);
        \draw[->>-=.5,black!50!green] (1.2,1) -- (0,1);
        \draw[->>-=.5,black!50!green] (-1.2,1) -- (0,1);
        \draw[->>-=.6,black!50!green] (0,-0.5) -- (1.2,-0.5);
        \draw[->>-=.6,black!50!green] (0,-0.5) -- (-1.2,-0.5);
        \filldraw[teal] (0,0) circle (2pt);
        \filldraw[black] (0,1) circle (2pt);
        \draw[] (-0.15,-0.15) node[] {O};
        \draw[] (0.15,1.15) node[] {Q};
			\end{tikzpicture}
   \caption{$\eps_2 \to \infty $.}
			\label{fig:K1_R1_rescaled_singular_eps_infty}
\end{subfigure}
			\caption{Limits of the singular dynamics of \eqref{eq:rob_reaction_2D_R2_rescaled}.}
			\label{fig:K1_R1_rescaled_singular}
    \end{figure}

The parameter $\eps_2$ changes the geometry of the critical manifold $\mathcal{S}$, compare
Figure \ref{fig:K1_R2_rescaled_singular} where $\mathcal{S}$ is shown in blue.
These changes are due to the occurrence of a transcritical bifurcation of $\mathcal{S}$ at $(y,z)^T=(-\sqrt{c},2c)^T$ for $\frac{1}{\eps_2^2}=c$.
We do not study this in detail since it occurs in the nonphysical part of phase space.
For  $\varepsilon_2 \in \Big[\sqrt{\frac{1}{\beta_2}}, \sqrt{\frac{1}{\beta_3}}\Big]$
these changes do not affect normal hyperbolicity of  $\mathcal{S}^a$. For $\eps_2 \to 0$, however, the fold point 
of  $\mathcal{S}$ approaches the equilibrium $Q$. In the limit $\eps_2=0$ the critical manifold is given by 
$$z=c-y^2,$$ 
i.e., the fold point of the critical manifold coincides with the equilibrium $Q=(0,c)^T$, see Figure \ref{fig:K1_R1_rescaled_singular_eps=0}. 
For $\eps_2 \to \infty$ the critical manifold $\mathcal{S}$
approaches the $y$-  and $z$-axis, see Figure \ref{fig:K1_R1_rescaled_singular_eps_infty}. In these two limits, normal hyperbolicity of $\mathcal{S}^a$ is lost at $Q$ and $O$, respectively.

Since we stay away from these degenerate limits in region $\Bar{B}_2$, the following construction of singular orbits 
and proof of their persistence based on  Fenichel theory works for all $\eps_2$ in region $\Bar{B}_2$.
In particular, the compact part of $\mathcal{S}^a$ connecting $y^{max}$ and the equilibrium $Q=(0,c)^T$ 
\begin{equation}
    \gamma_0^s= \{(y,z)^T \in \mathcal{S}^a:y \in [0,y^{max}]\}.
\end{equation} 
is normally attracting (which is indicated by green double arrows in Figure \ref{fig:K1_R2_rescaled_singular}). The first part of the singular orbit, which connects the initial value $O=(0,0)^T$ along the fast fiber (green) with the point $(y^{max},0)^T \in \mathcal{S}^a$ is given by 
\begin{equation}
    \gamma_0^f=\{(y,0)^T \in \R^2:y \in [0,y^{max}]\} 
\end{equation} 
It remains to check the reduced flow on $\mathcal{S}^a$. We change to the slow time scale $t=r \tau$ and obtain the reduced flow on $\mathcal{S}^a$
\begin{equation} \label{eq:K12_slow-flow}
    \Dot{z} = y^2 \geq 0,
\end{equation}
where \enquote{ $\Dot{ }$ } denotes differentiation w.r.t. $t$. Thus, the solution of the reduced problem starting at $(y^{max},0)^T$ converges center-like, i.e., with an algebraic rate, to the equilibrium $Q=(0,c)^T$. We obtain the following lemma.
\begin{lem} \label{lem:singular orbit B2}
    There exists a singular orbit $\gamma_0:=\gamma_0^f\cup \gamma_0^s$ of \eqref{eq:rob_reaction_2D_R2_rescaled_singular} connecting the initial value $O$ and the equilibrium $Q$.
\end{lem}
Due to the following theorem, the singular orbit perturbs to a genuine orbit for $r$ small.
\begin{thm} \label{thm:convergence_R_2}
    There exists $r_0>0$ such that for all $\varepsilon_2 \in \Big[\sqrt{\frac{1}{\beta_2}}, \sqrt{\frac{1}{\beta_3}}\Big]$ and $r \in (0,r_0]$ there exists a smooth orbit $\gamma_r$ of system \eqref{eq:rob_reaction_2D_R2_rescaled}, connecting the initial value $O=(0,0)^T$ with the equilibrium $Q=(0,c)^T$. The perturbed orbit $\gamma_r$ is $\mathcal{O}(r)$-close to $\gamma_0$ in Hausdorff distance.
\end{thm} 
\begin{proof}
The normally hyperbolic attracting  critical manifold $\mathcal{S}^a$ perturbs to an attracting slow manifold $\mathcal{S}^a_r$ by Fenichel theory for $0 < r \ll 1$, which contains the equilibrium $Q$. Since there are no further
equilibria for $r >0$ the slow flow on $\mathcal{S}^a_r$ converges to $Q$ for $y \geq 0$ (as a center flow).
Viewed as an equilibrium of  system \eqref{eq:rob_reaction_2D_R2_rescaled} $Q$ has a two-dimensional center-stable manifold
$W^{cs}$ which intersects $\mathcal{S}^a_r$ transversally. 
By Fenichel theory the solution with inital value $O$, i.e. the orbit $\gamma_r$,  is attracted exponentially 
onto $\mathcal{S}^a_r$ and hence converges to $Q$ in the center direction. By construction $\gamma_r$ is $\mathcal{O}(r)$ close to $\gamma_0$.
\end{proof}

We conclude that for all $(\eps_1,\eps_2)^T \in B_2$ with $\vert \vert (\eps_1,\eps_2)^T \vert \vert <r_0$ there exists a smooth orbit $\gamma_{\eps_1}$ of system \eqref{equ:rob_reaction_2D}, which is $\mathcal{O}(\sqrt{\eps_1})$-close to $\gamma_0$ in Hausdorff distance, connecting the initial value $O$ with the equilibrium $Q$.

A possibility to compare our asymptotic results with the numerics is the maximal value of the $y$-component $y^{max}$, which we will focus on in the following. Due to the extra rescaling \eqref{eq:rescaling_R2} we even achieve an error estimate of $\mathcal{O}(\eps_1)$ in $y$-direction. Indeed, by undoing the rescalings \eqref{equ:parameter_blow-up-chart_K1} and \eqref{eq:rescaling_R2} we obtain
\begin{equation} \label{eq:Evaluate_numerics_K_2}
    y=r\Tilde{y}=\sqrt{\eps_1}\Tilde{y}.
\end{equation}
Along the singular orbit $\gamma_0$ we have 
$$\max \{\Tilde{y}:(\Tilde{y},z)^T  \in \gamma_0\}=\Tilde{y}^{max}=\sqrt{c},$$
such that
$$y^{max}=\sqrt{\eps_1}(\sqrt{c}+\mathcal{O}(\sqrt{\eps_1}))=\sqrt{\eps_1 c}+\mathcal{O}(\eps_1).$$
Inserting the parameter values \eqref{equ:rob_rates} of the original problem gives
\begin{equation} \label{eq:y_max}
    y^{max}=3.651\cdot 10^{-5}+\mathcal{O}(10^{-9}),
\end{equation}
which fits well with the value obtained by numerical simulations, e.g., compare with Figure \ref{fig:rob_numeric}. In particular, this confirms the numerical results in \cite{Hairer_1996}. 

\comm{\begin{rem}
From the blow-up point of view the rescaling \eqref{eq:rescaling_R2} can be viewed as the scaling chart of a cylindrical blow-up
of the degenerate line $(0,z,0)$, $z \in \R$ in extended $(y,z,r)$ phase space. Since this scaling chart covers the relevant dynamics, there is no need to introduce this blow-up explicitly.
\end{rem}}

It remains to do the analysis of the degenerate cases corresponding to $\Tilde{\eps}_2 \to 0$ and $\Tilde{\eps}_2 \to \infty$ in regions $\Bar{B}_1$ and $\Bar{B}_3$, respectively. We start with the region $\Bar{B}_1$, since this allows us to continue in the current chart $\mathcal{P}_1$.

\section{Analysis in Region $B_1$} \label{sec:Regime R1}
The starting point of the following analysis in chart $\mathcal{P}_1$ are the equations \eqref{eq:rob_reaction_2D_R2_rescaled}, which we restate here for convenience
\begin{equation} \label{eq:Bar{B}_1_rescaled}
    \begin{aligned}
y'&=c-r y-z - y^2 -  \Tilde{\eps}_2 yz\\
z'&= r y^2.
\end{aligned}
\end{equation}
As described before, for $\Tilde{\eps}_2 \to 0$ the fold point of $\mathcal{S}$ is at $Q=(0,c)^T$, see Figure \ref{fig:K1_R1_rescaled_singular_eps=0}. To treat this loss of normal hyperbolicity we perform a blow-up of the fold point $(0,c,0)^T$ in extended $(y,z,\Tilde{\eps}_2)^T$ space. To handle the terms $-ry$ and $-\Tilde{\eps}_2yz$ in \eqref{eq:Bar{B}_1_rescaled}, an additional homogeneous parameter blow-up of the origin $(r,\Tilde{\eps}_2)^T=(0,0)^T$ in chart $\mathcal{P}_1$ is introduced. Otherwise, we would not be able to desingularize the dynamics in the blow-up of the fold point. In the original parameters the second parameter blow-up amounts to dividing the region $B_1$ into two parts $B_{11}$ and $B_{12}$ by a curve 
\begin{equation} \label{eq:curve_C1}
    C_1:=\{(\eps_1,\eps_2)^T \in \R^2:\eps_1=\beta_1 \eps_2\}
\end{equation}
with some $\beta_1>0$. The parts $B_{11}$ and $B_{12}$ correspond to scaling regimes $0\leq \eps_2 \lesssim \eps_1$ and $\eps_2^2 \ll \eps_1 \lesssim \eps_2$, respectively. 
The second parameter blow-up map is given by 
\begin{equation} \label{equ:parameter_blow-up_2}
 \begin{aligned}
         \Phi_{par}^2: [0,\infty) \times \mathbb{S}^1 &\to \R^2\\
         (s,\Bar{r},\Bar{\Tilde{\eps}}_2) &\mapsto \begin{cases}
             r=s \Bar{r}\\
             \Tilde{\eps}_2=s \Bar{\Tilde{\eps}}_2.
         \end{cases}
 \end{aligned}
\end{equation}

In the blown-up space the curve $C_1$ corresponds to the line
$$\Bar{\eps}_2=\sqrt{\frac{1}{1+\beta_1^2}}.$$

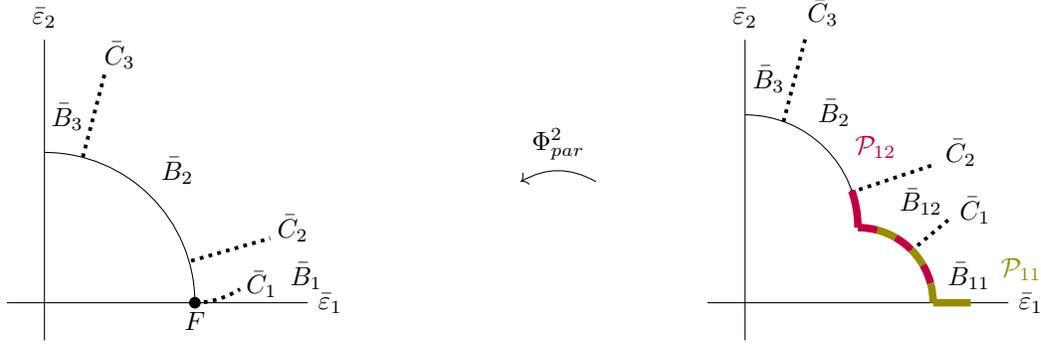
\begin{figure}
\centering
\begin{subfigure}{0.4\textwidth}
    \begin{tikzpicture}[scale=1]
        \draw[] (0,-0.5) -- (0,3.5);
        \draw[] (-0.5,0) -- (3.5,0);
        \draw[] (2, 0) arc (0:90:2);
        \draw[scale=1, domain=1.95:3, smooth, variable=\x, line width=0.5mm, dotted] plot ({\x}, {(1/3.5)*\x});
        \draw[scale=1, domain=0.51:0.8, smooth, variable=\x, line width=0.5mm, dotted] plot ({\x}, {(3.8)*\x});
        \draw[scale=1, domain=2:2.6, smooth, variable=\x, line width=0.5mm, dotted] plot ({\x}, {(1/2)*(\x-2)^2});
        \draw (83:2.5) node{$\Bar{B}_3$};
        \draw (45:2.5) node{$\Bar{B}_2$};
        \draw (6:3.5) node{$\Bar{B}_1$};
        \filldraw[] (2,0) circle (2pt);
        \draw (3.8,0) node{$\Bar{\eps}_1$};
        \draw (0,3.8) node{$\Bar{\eps}_2$};
        \draw[] (5:2.9) node{$\Bar{C}_1$};
        \draw[] (3.3,1) node{$\Bar{C}_2$};
        \draw[] (1,3.3) node{$\Bar{C}_3$};
        \draw[below] (2,0) node{$F$};

			\end{tikzpicture}
\end{subfigure}
\begin{subfigure}{0.15\textwidth}
\begin{tikzpicture}
       \path[->]
    (1,2) edge[bend right] node [left] {} (0,2);
    \draw (0.5,2.5) node{$\Phi_{par}^2$};
    \draw (0,0) node{};
\end{tikzpicture}
\end{subfigure}
\begin{subfigure}{0.4\textwidth}
\begin{tikzpicture}[scale=1]
        \draw[] (0,-0.5) -- (0,3.5);
        \draw[] (-0.5,0) -- (3.5,0);
        \draw[] (0, 2.5) arc (90:0:1.5);
        \draw[] (2.5,0) arc (0:90:1);
        \draw[scale=1, domain=1.4:2.5, smooth, variable=\x, line width=0.5mm, dotted] plot ({\x}, {(1/3)*\x+1});
        \draw[scale=1, domain=0.51:0.8, smooth, variable=\x, line width=0.5mm, dotted] plot ({\x}, {(3.8)*\x+0.47});
        \draw[scale=1, domain=2.2:2.7, smooth, variable=\x, line width=0.5mm, dotted] plot ({\x}, {(1/1.2)*\x-1.15});
        \draw (84:3) node{$\Bar{B}_3$};
        \draw (65:2.8) node{$\Bar{B}_2$};
        \draw (30:2.7) node{$\Bar{B}_{12}$};
        \draw (7:3) node{$\Bar{B}_{11}$};
        \draw[olive,line width=1mm] (2.5, 0) arc (0:15:1);
        \draw[olive,line width=1mm] (2.5,0) -- (3,0);
        \draw[purple,line width=1mm] (1.5,1) arc (90:75:1);
        \draw[purple,domain=0:19,line width=1mm] plot ({1.5*cos(\x)}, {1.5*sin(\x)+1});
        \draw[olive,domain=60:75,line width=1mm] plot ({cos(\x)+1.5}, {sin(\x)});
        \draw[olive,domain=30:45,line width=1mm] plot ({cos(\x)+1.5}, {sin(\x)});
        \draw[purple,domain=60:45,line width=1mm] plot ({cos(\x)+1.5}, {sin(\x)});
        \draw[purple,domain=30:15,line width=1mm] plot ({cos(\x)+1.5}, {sin(\x)});
        \draw[purple,line width=1mm] (50:2.7) node{$\mathcal{P}_{12}$};
        \draw[olive,line width=1mm] (7:3.7) node{$\mathcal{P}_{11}$};
        \draw[] (22:3.3) node{$\Bar{C}_1$};
        \draw[] (35:3.5) node{$\Bar{C}_2$};
        \draw[] (75:4) node{$\Bar{C}_3$};
        \draw (3.8,0) node{$\Bar{\eps}_1$};
        \draw (0,3.8) node{$\Bar{\eps}_2$};

			\end{tikzpicture}
\end{subfigure}
			\captionof{figure}{Schematic picture of the parameter blow-up $\Phi_{par}^2$ of the point $F$ (corresponding to the origin in chart $\mathcal{P}_1)$ shown in blown-up $(\R \times \mathbb{S}^1)$-space.}
			\label{fig:second_parameter_blow-up_P1}
\end{figure}

Again it is convenient to perform the analysis in two charts corresponding to the directions $\Bar{r}=1$ and $\Bar{\eps}_2=1$, respectively. In these charts the blow-up transformation $\Phi_{par}^2$ is given by 
\begin{align} 
         \mathcal{P}_{11}&: r=s,\, \tilde{\eps}_2=s \eps_{21}, \label{eq:parameter_blow-up_2_K11}\\
         \mathcal{P}_{12}&: r=sr_1, \tilde{\eps}_2=s. \label{eq:parameter_blow-up_2_K13}
 \end{align}
 
such that in chart $\mathcal{P}_{11}$ the regions $B_{11}$ and $B_{12}$ in the blown-up space are given by $\eps_{21}<\beta_1$ and $\eps_{21}\geq \beta_1$, respectively.
A schematic representation of the second blow-up in parameter space is shown in Figure \ref{fig:second_parameter_blow-up_P1}. As can be seen in Figure \ref{fig:second_parameter_blow-up_P1}, chart $\mathcal{P}_{11}$ will be used for analysing the limit $\tilde{\eps}_2 \to 0$ in region $\Bar{B}_{11}$, whereas chart $\mathcal{P}_{12}$ covers region $\Bar{B}_{12}$. We begin with the analysis in chart $\mathcal{P}_{11}$. 

\subsection{Analysis in region $B_{11}$}

Inserting the parameter blow-up transformation \eqref{eq:parameter_blow-up_2_K11} into \eqref{eq:Bar{B}_1_rescaled} we obtain 
\begin{equation} \label{eq:Bar{B}_1_rescaled_K11}
    \begin{aligned}
y'&=c-s y-z - y^2 -  s\eps_{21} yz\\
z'&= s y^2.
\end{aligned}
\end{equation}

In the following analysis, it will be important that $\eps_{21} \in [0,\beta_1]$. System \eqref{eq:Bar{B}_1_rescaled_K11} is of classical slow-fast type with parameter $s$. The critical manifold is given by  
$$\mathcal{S}=\{(y,z)^T \in \R^2: z=c-y^2\}$$
with a fold point at the equilibrium $Q=(0,c)^T$, see again Figure \ref{fig:K1_R1_rescaled_singular}. The candidate singular orbit starting from the initial value $O=(0,0)^T$ is again $\gamma_0=\gamma_0^f \cup \gamma_0^s$, but it approaches $Q$ along the slow manifold $\mathcal{S}$, which loses normal hyperbolicity at the fold point. Therefore, we cannot use Fenichel theory directly to prove convergence to the genuine equilibrium along $\gamma_0$. We resolve this degeneracy by artificially adding $s'=0$ to \eqref{eq:Bar{B}_1_rescaled_K11} and applying a spherical blow-up of the nilpotent point $(y,z,s)^T=(0,c,0)^T$ of this extended system, see \cite{Krupa_2001_Extend} for a detailed explanation of the blow-up method in the context of planar fold points.

The suitable blow-up transformation is given by 
\begin{equation} \label{eq:spherical_blow-up_R1}
     \begin{aligned}
         \Phi:[0,\infty) \times \mathbb{S}^2  & \to \R^3\\
	(\sigma,\bar{y},\bar{z},\bar{s}) & \mapsto \begin{cases}
	    y=\sigma \bar{y}\\
        z=c+\sigma^2\bar{z}\\
        s=\sigma \bar{s}.
	\end{cases}
     \end{aligned}
 \end{equation}

\begin{figure}
    \centering
    \includegraphics{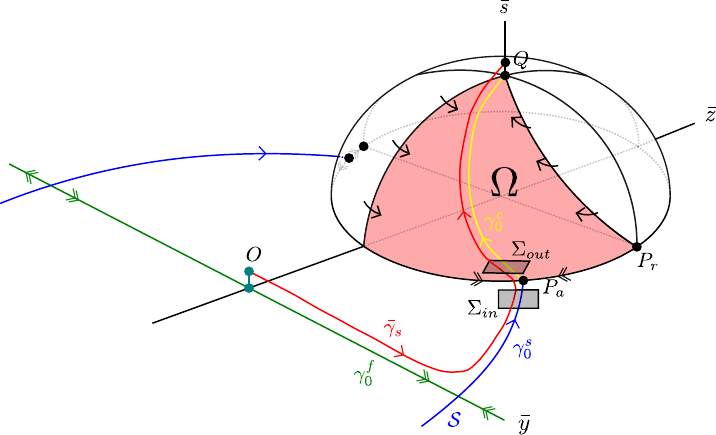}
    \caption{Dynamics of the blown-up extended system \eqref{eq:Bar{B}_1_rescaled_K11}.}
    \label{fig:K_11_singular dynamics}
\end{figure}

The nilpotent point $(0,c,0)^T$ is blown-up to the sphere $\{0\} \times \mathbb{S}^2$, which is the preimage of $(0,c,0)^T$ under the map $\Phi$, see Figure \ref{fig:K_11_singular dynamics}.

 \begin{rem} \label{rem:P11_different_weights_than_fold}
     Note that in the transformation \eqref{eq:spherical_blow-up_R1} the weights of the radial variable $\sigma$ deviate from the weights in the analysis of the generic fold point. This is a consequence of the fold point coinciding with an equilibrium in our case.
 \end{rem}

Much of the following analysis proceeds along the lines of \cite{Krupa_2001_Extend}, the dynamics on the sphere $\sigma=0$ is, however, different from the standard fold point, so we give the necessary details.
Again it will be convenient to work in directional charts which correspond to directions $\Bar{y}=1$, $\Bar{s}=1$, and $\Bar{z}=-1$. The blow-up transformation in these charts is given by 
   \begin{align} 
        &\mathcal{K}_{11}^1: y=\sigma_1, \quad z=c+\sigma_1^2z_1, \quad s=\sigma_1s_1\label{eq:spherical_blow-up_K^1}\\
        &\mathcal{K}_{11}^2: y=\sigma_2y_2, \quad z=c+\sigma_2^2z_2, \quad s=\sigma_2 \label{eq:spherical_blow-up_K^2}\\
        &\mathcal{K}_{11}^3: y=\sigma_3y_3, \quad z=c-\sigma_3^2, \quad s=\sigma_3s_3,\label{eq:spherical_blow-up_K^3}
    \end{align}
respectively. Note that subscripts refer to the parameter blow-up chart, whereas a superscript denotes the corresponding chart in phase space. Chart $\mathcal{K}_{11}^1$ covers the right ($\Bar{y}>0$) side of the sphere, chart $\mathcal{K}_{11}^2$ covers the top ($\Bar{s}>0$) of the sphere, and chart $\mathcal{K}_{11}^3$ covers the front ($\Bar{z}<0$) side of the sphere, see Figure \ref{fig:K_11_singular dynamics}.

\begin{rem} \label{chartnotation}
For blow-ups in phase space we will often follow the useful convention, that an object $A$ is denoted as
$A_i$ in a chart $\mathcal{K}^i$, $i=1,2,3$ in which the blow-up is studied. As an example consider the equilibrium $Q$ which will be studied in chart $\mathcal{K}_{11}^2$ and is denoted as $Q_2$ there.
\end{rem}

The following subset of the sphere is central for our analysis.

\begin{defin} \label{def:omega}
    Let $\Omega$ be the compact subset of the sphere $(\sigma=0)$ enclosed by the equator $(\Bar{s}=0)$, the meridian $(\Bar{y}=0)$, and the curve which is represented by $s_1=-z_1$ in chart $\mathcal{K}_{11}^1$ and $z_2=-y_2$ in chart $\mathcal{K}_{11}^2$, see Figure \ref{fig:K_11_singular dynamics} where $\Omega$ is shown in red.
\end{defin}

We have the following result.

\begin{lem} \label{rem:plan_phase_plane_argument}
    The flow of the blown-up vector field on the sphere has the properties:
    \begin{enumerate}
        \item The set  $\Omega$ is forward invariant.
        \item There exists a heteroclinic orbit $\gamma_0^c$ connecting the endpoint $P_a$ of $\mathcal{S}$ with the equilibrium $Q$.
    \end{enumerate}

\end{lem}

\begin{proof}
We start the analysis in chart $\mathcal{K}_{11}^1$, which is one of the entrance charts since it contains the endpoint $P_a$ of the attracting branch of the critical manifold $\mathcal{S}$, denoted by $\mathcal{S}^a$, with reduced flow towards the sphere. Inserting \eqref{eq:spherical_blow-up_K^1} into \eqref{eq:Bar{B}_1_rescaled_K11} and after desingularizing, i.e., dividing out a factor of $\sigma_1$, we obtain
\begin{equation} \label{eq:K_11_blow-up_K^1_dynamics}
    \begin{aligned}
z_1'&= s_1+2z_1(s_1+z_1+1+\sigma_1^2\eps_{21}s_1z_1+s_1\eps_{21}c)\\
s_1'&=s_1(s_1+z_1+1+\sigma_1^2\eps_{21}s_1z_1+s_1\eps_{21}c)\\
\sigma_1'&=-\sigma_1(s_1+z_1+1+\sigma_1^2\eps_{21}s_1z_1+s_1\eps_{21}c).
\end{aligned}
\end{equation}
The planes $\sigma_1=0$ and $s_1=0$ are invariant. They intersect in a line, which corresponds to a part of the equator of the sphere, on which the dynamics is governed by $z_1'=2z_1(z_1+1)$. There are two equilibria $P_a=(-1,0,0)^T$ and $P_r=(0,0,0)^T$ which are attracting and repelling on this line with eigenvalues $-2$ and $2$, respectively.

On the plane $s_1=0$ the dynamics is given by 
\begin{equation} \label{eq:K_11_blow-up_K^1_s=0}
    \begin{aligned}
    z_1'&= 2z_1(z_1+1)\\
    \sigma_1'&=-\sigma_1(z_1+1).
\end{aligned}
\end{equation}
The normally attracting line of equilibria 
$$z_1=-1$$
corresponds to the attracting branch of the critical manifold $\mathcal{S}^a$, see Figure \ref{fig:K_11_singular dynamics}.  We now investigate the dynamics on the plane $\sigma_1=0$ (on the sphere) near the point $P_a$ governed by 
\begin{equation} \label{eq:K_11_blow-up_K^1_sigma=0}
    \begin{aligned}
    z_1'&= s_1+2z_1(s_1+z_1+1+s_1\eps_{21}c)\\
s_1'&=s_1(s_1+z_1+1+s_1\eps_{21}c).
\end{aligned}
\end{equation}

We recover the equilibria $P_a$ and $P_r$. The eigenvalues of the linearization at $P_a$ and $P_r$ are $-2$, $0$ and $2$, $1$, respectively. We conclude that $P_r$ is a source on the sphere. Standard center manifold theory \cite{Guckenheimer_1983} implies the existence of an attracting one-dimensional center manifold $N_a$ at $P_a$, which is given as a graph $z_1=h_1(s_1)$ with expansion
\begin{equation} \label{eq:K_11^1_center_expansion_N_a}
    h_1(s_1)=-1-(\frac{1}{2}+\eps_{21}c)s_1+ \mathcal{O}(s_1^2).
\end{equation}
The corresponding flow on $N_a$ is governed by $$s_1'=s_1^2/2+\mathcal{O}(s_1^3),$$ 
hence $z_1$ increases along $N_a$. This implies that the branch of $N_a$ in $s_1>0$ is unique. 
For proving assertion (ii), it remains to show that the continuation of this branch of $N_a$ by the flow connects $P_a$ with the equilibrium $Q$ (which is only visible in the scaling chart $\mathcal{K}_{11}^2$). This is done in the following by first proving assertion (i) and using a phase plane argument. 

The part of $\partial \Omega$ that is visible in chart $\mathcal{K}_{11}^1$ is given by the invariant half line $s_1=0$ and the line $s_1=-z_1$, respectively, both with $z_1\leq 0$. On these half lines the flow cannot exit $\Omega$. For the line $s_1=-z_1$ this follows from 
\begin{equation} \label{eq:K_11^2_trapping_2D}
    (s_1+z_1)'|_{s_1=-z_1}=-s_1^2\eps_{21}c \leq 0.
\end{equation}
for all $\eps_{21}\geq 0$, see \cite[p.~219]{Amann_1990}.

Now we switch to the chart $\mathcal{K}_{11}^3$ in which the governing equations are
\begin{equation} \label{eq:K_11_blow-up_K^3_dynamics}
    \begin{aligned}
y_3'&=-y_3s_3+1-y_3^2+\sigma_3^2\eps_{21}y_3s_3-s_3\eps_{21}y_3c+\frac{1}{2}y_3^3s_3\\
s_3'&= \frac{1}{2}y_3^2s_3^2\\
\sigma_3'&=-\frac{1}{2}\sigma_3y_3^2s_3.
\end{aligned}
\end{equation}

On the invariant plane $s_3=0$ we recover the two normally hyperbolic parts of the critical manifold as lines of equilibria 
$$y_3=\pm 1,$$
where the attracting line $y_3=1$ terminates in $P_a$. Clearly this chart also covers the center manifold $N_a$ originating at $P_a$.

The part of $\partial \Omega$ on the sphere $\sigma=0$ that is visible in chart $\mathcal{K}_{11}^3$ is given by the invariant half line $s_3=0$, $y_3\geq 0$ and the half line $y_3=0$, $s_3\geq 0$. The flow on the sphere cannot leave $\Omega$ at these half lines. For the line $y_3=0$, $s_3\geq 0$ this follows from $y_3'=1$.

To cover the part of $\Omega$ close to $Q$ we change to the scaling chart $\mathcal{K}_{11}^2$ where we can trace $\gamma_0^c$ once it has entered $\Omega$. The dynamics in the scaling chart $\mathcal{K}_{11}^2$ is governed by
\begin{equation} \label{eq:K_11_blow-up_K^2_dynamics}
    \begin{aligned}
y_2'&=-y_2-z_2 - y_2^2 -  y_2\eps_{21}c-\sigma_2^2z_2y_2\eps_{21}\\
z_2'&= y_2^2\\
\sigma_2'&=0.
\end{aligned}
\end{equation}

On the invariant sphere $\sigma=0$ this simplifies to
\begin{equation} \label{eq:K_11_blow-up_K^2_dynamics_singular_2D}
    \begin{aligned}
y_2'&=-y_2-z_2 - y_2^2 -  y_2\eps_{21}c\\
z_2'&= y_2^2.
\end{aligned}
\end{equation}

The boundary of $\Omega$ in $\Bar{s}>0$ is given by parts of the lines $y_2=0$, $z_2\leq 0$ and $z_2=-y_2$, $y_2\geq 0$. The flow of \eqref{eq:K_11_blow-up_K^2_dynamics_singular_2D} cannot leave $\Omega$ at these parts of the boundary since $$y_2'=-z_2$$
on the line $y_2=0$ and 
$$(y_2+z_2)'=-y_2\eps_{21}c\leq 0$$
on the line $z_2=-y_2$. We conclude that $\Omega$ is indeed a compact forward invariant trapping region on the sphere. This concludes the proof of assertion (i). 

In chart $\mathcal{K}_{11}^2$ the equilibrium $Q$ corresponds to the point $Q_2=(0,0)^T$. The linearization of \eqref{eq:K_11_blow-up_K^2_dynamics_singular_2D} at $Q_2$ has eigenvalues $\lambda_s=-1-\eps_{21}c$ and $\lambda_c=0$ with corresponding eigenvectors $v_s=(1,0)^T$ and $v_c=(1,-1-\eps_{21}c)^T$. Standard center manifold theory implies the the existence of an attracting (non-unique) center manifold $W^c$, which lies in the interior of $\Omega$ for $\eps_{21}>0$ and coincides with the line $z_2=-y_2$ in the limiting case $\eps_{21}=0$. The flow on the center manifold $W^c$ in $\Omega$ is directed towards the equilibrium $Q_2$. There are no equilibria in the interior of $\Omega$ such that we can exclude periodic orbits. The only equilibrium in the forward invariant compact set $\Omega$ which is not repelling is the equilibrium $Q_2 \subset \partial \Omega$. On the sphere $\sigma=0$, the Poincare-Bendixson theorem applies, therefore all orbits within $\Omega$ must converge to the equilibrium $Q_2$ along the center manifold $W^c$. Therefore, the continuation of the $N_a$ in $s_1 >0$ converges to $Q_2$.
We denote the corresponding heteroclinic orbit by $\gamma_0^c$, which is shown in yellow in Figure \ref{fig:K_11_singular dynamics}. This proves assertion (ii).
\end{proof}

By collecting the results of this subsection we obtain.

\begin{lem} \label{lem:singular_orbit_P11}
    There exists a singular orbit $\gamma_0$ of the blown-up extended system \eqref{eq:Bar{B}_1_rescaled_K11} connecting the initial value $O$, via $P_a$, with the equilibrium $Q$. 
\end{lem}
\begin{proof}
    Starting from the initial value $O$, we follow $\gamma_0^f$ and $\gamma_0^s$ as before. In the blown-up problem $\gamma_0^s$ terminates in the point $P_a$. From there we follow $\gamma_0^c$ which connects $P_a$ and $Q$.
    We define the singular orbit as $\gamma_0:=\gamma_0^f \cup \gamma_0^s \cup \gamma_0^c$, see Figure \ref{fig:K_11_singular dynamics}.
\end{proof}

Now we prove that the singular orbit $\gamma_0$ perturbes to a smooth orbit $\gamma_s$ connecting the initial condition $O$ and the equilibrium $Q$ for $0<s \ll 1$.

\begin{thm} \label{thm:P11}
    There exists a constant $\Tilde{s}>0$ such that for all $\eps_{21} \in [0,\beta_1]$ and $s \in (0,\Tilde{s}]$ , there exists a smooth orbit $\gamma_{s}$ of \eqref{eq:Bar{B}_1_rescaled_K11} connecting the initial value $O$ and the equilibrium $Q$. The corresponding orbit $\Bar{\gamma}_s$
    in blown-up space  is $\mathcal{O}(s)$-close to $\gamma_0$ in Hausdorff distance. 
\end{thm}
\begin{proof}
    The proof is carried out in the blow-up of system \eqref{eq:Bar{B}_1_rescaled_K11} extended by $s'=0$. In a first step we show that the continuation of the slow manifold by the flow, which exist by Fenichel theory away from the fold, converges to the equilibrium $Q$. For this purpose we define two sections in the entrance chart $\mathcal{K}_{11}^1$ close to $P_a$ as
\begin{equation} \label{eq:incoming_section}
    \Sigma_{in}:=\{(z_1,s_1,\sigma_1)^T \in \R^3: \sigma_1=a,\, |z_1+1|<b,\, s_1<a\}
\end{equation}
 and
\begin{equation} \label{eq:outgoing_section}
    \Sigma_{out}:=\{(z_1,s_1,\sigma_1)^T \in \R^3: \sigma_1<a,\, |z_1+1|<b,\, s_1=a\}
\end{equation}
with $a,b>0$ small enough, see Figure \ref{fig:K_11_singular dynamics}. Away from the sphere $\sigma=0$ the attracting branch $\mathcal{S}^a$ of the critical manifold perturbes to an attracting slow manifold $\mathcal{S}^a_s$ for $s \ll 1$ by Fenichel theory. In extended phase space, this one-parameter family of slow manifolds can be viewed as a two-dimensional invariant attracting slow manifold $\mathcal{M}$. The manifold $\mathcal{M}$ is defined at least up to the section $\Sigma_{in}$. We extend the manifold $\mathcal{M}$ by the forward flow of the blown-up vector field past $\Sigma_{in}$, the corresponding larger manifold is still denoted as $\mathcal{M}$. The results in \cite{Krupa_2001_Extend} on the standard singularly perturbed fold point imply that $\mathcal{M}$ is attached to the orbit $\gamma_0^c$. Therefore, we can track $\mathcal{M}$ across the sphere for $s$ small. In the blown-up phase space the equilibrium $Q$ corresponds to a line of equilibria $(0,0,\sigma_2)^T$, $\sigma_2 \in [0,\Tilde{s}]$. The linearization along this line of equilibria has one negative and a double zero eigenvalue. Standard invariant manifold theory implies the existence of a three-dimensional center-stable manifold $W^{cs}$ of this line of equilibria. Since the orbit $\gamma_0^c$ on the sphere intersects $W^{cs}$ transversely, the manifold $\mathcal{M}$ also intersects $W^{cs}$ since it is a small smooth perturbation of $\gamma_0^c$ for $s \ll 1$. This implies that all orbits in $\mathcal{M}$ converge to $Q$. 

Viewed in chart $\mathcal{K}_{11}^3$ of the extended blown-up phase space, the line of initial conditions $\{ (0,0,s)^T,  \; s \in [0,\Tilde{s}] \} $  corresponds to the line $(0,\sqrt{c},\frac{s}{\sqrt{c}})^T$, $s \in [0,\Tilde{s}]$. All orbits starting on this line are exponentially attracted onto the manifold $\mathcal{M}$ by Fenichel theory until they reach $\Sigma_{in}$. During the passage from $\Sigma_{in}$ to $\Sigma_{out}$ an additional exponential contraction towards $\mathcal{M}$ occurs due to \cite[Proposition 2.8]{Krupa_2001_Extend}. Beyond the section $\Sigma_{out}$ the evolution of these orbits is governed by system \eqref{eq:K_11_blow-up_K^2_dynamics} with $\sigma_2 \in [0,\Tilde{s}]$. Since $\sigma_2$ acts as a regular perturbation parameter, these orbits intersect $W^{cs}$ for $\Tilde{s}$ sufficiently small. This implies the existence of a smooth perturbed orbit $\Bar{\gamma}_s$ connecting the lines of equilibria corresponding to $O$ and $Q$, respectively. The assertions of the theorem follow by applying the blow-up transformation \eqref{eq:spherical_blow-up_R1}, i.e., 
$\gamma_s = \Phi( \Bar{\gamma}_s)$.

\end{proof}

In order to complete the argument in region $\Bar{B}_1$, we use chart $\mathcal{P}_{12}$ which covers the region $\Bar{B}_{12}$.

\subsection{Analysis in region $B_{12}$}

We insert the transformation \eqref{eq:parameter_blow-up_2_K13} into \eqref{eq:Bar{B}_1_rescaled} and obtain 
\begin{equation} \label{eq:Bar{B}_1_rescaled_K13}
    \begin{aligned}
y'&=c-s r_1y-z - y^2 -  s yz\\
z'&= s r_1y^2.
\end{aligned}
\end{equation}
where $s$ and $r_1$ are both small. The goal is to construct the orbit connecting $O$ to the equilibrium $Q$ for $(s,r_1)^T$, $r_1 >0$, $s>0 $ in a small neighborhood of the origin.
For $s=0$ we again obtain the critical manifold 
$$z=c-y^2$$
with  the equilibrium $Q=(0,c)^T$ at the fold point.

Again we use the blow-up transformation \eqref{eq:spherical_blow-up_R1} to resolve the degeneracy of the fold point and we obtain the following result.

\begin{lem} \label{lem:S^a_B_12}
    In the blown-up space of  system \eqref{eq:Bar{B}_1_rescaled_K13} extended by the equation $s'=0$ there exists for $r_1 =0$ a two-dimensional attracting critical manifold $\mathcal{S}^a$ (blue in Figure \ref{fig:K_13_singular dynamics}), which contains the line of equilibria corresponding to the genuine equilibrium $Q$. The critical manifold $\mathcal{S}^a$ perturbes regularly to a slow manifold $\mathcal{S}^a_{r_1}$ for $r_1$ small enough. All orbits of the reduced flow on $\mathcal{S}^a$ approach the line of equilibria in the center direction, i.e. with an algebraic rate,  for all $0<s \leq 1/\sqrt{\beta_2}$.
\end{lem}

\begin{proof}
    We carry out the analysis in two directional charts
\begin{align} 
        &\mathcal{K}_{12}^2: y=\sigma_2y_2, \quad z=c+\sigma_2^2z_2, \quad s=\sigma_2 \label{eq:spherical_blow-up_K_13^2}\\
        &\mathcal{K}_{12}^3: y=\sigma_3y_3, \quad z=c-\sigma_3^2, \quad s=\sigma_3s_3\label{eq:spherical_blow-up_K_13^3}
    \end{align}
covering  the top $\Bar s >0$ and the front $\Bar{z} < c$ part of the sphere, respectively,
see  Figure \ref{fig:K_13_singular dynamics}. The parts of $\mathcal{S}^a$ investigated in chart $\mathcal{K}_{12}^2$ and $\mathcal{K}_{12}^3$ are denoted by $\mathcal{S}^a_2$ and $\mathcal{S}^a_3$, respectively.

\begin{figure}
    \centering
    \includegraphics{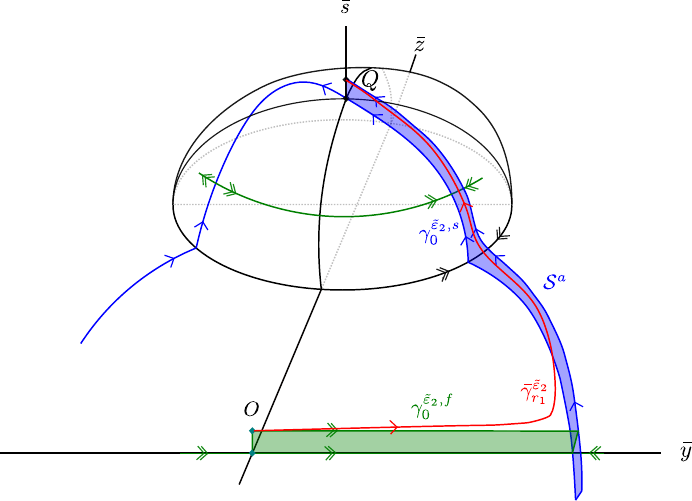}
    \caption{Dynamics of the blown-up extended system \eqref{eq:Bar{B}_1_rescaled_K13}.}
    \label{fig:K_13_singular dynamics}
\end{figure} 

As before we start the analysis in the entrance chart $\mathcal{K}_{12}^3$. By inserting \eqref{eq:spherical_blow-up_K_13^3} into \eqref{eq:Bar{B}_1_rescaled_K13} and desingularizing by dividing out a factor $\sigma_3$
we obtain
\begin{equation} \label{eq:K_13_blow-up_K^3_dynamics}
    \begin{aligned}
    y_3'&= 1-y_3^2-y_3s_3c+\sigma_3^2s_3y_3+r_1(s_3y_3+\frac{1}{2}s_3y_3^3)\\
    s_3'&=\frac{1}{2}r_1s_3^2y_3^2\\
    \sigma_3'&=-\frac{1}{2}r_1\sigma_3s_3y_3^2.
\end{aligned}
\end{equation}
Note that system \eqref{eq:K_13_blow-up_K^3_dynamics} is of standard slow-fast type with singular perturbation parameter $r_1$ and corresponding layer problem

\begin{equation} \label{eq:K_13_blow-up_K^3_r_1=0}
    \begin{aligned}
    y_3'&= 1-y_3^2-y_3s_3c+\sigma_3^2s_3y_3\\
    s_3'&=0\\
    \sigma_3'&=0.
\end{aligned}
\end{equation}

For $r_1=0$ we find the two-dimensional critical manifold
\begin{equation} \label{eq:crit_mfld_K13^3}
    \mathcal{S}_3=\{(y_3,s_3,\sigma_3)^T \in \R^3:1-y_3^2-y_3s_3c+\sigma_3^2s_3y_3=0 \},
\end{equation}
which for $s_3=0$ reduces to the two lines of equilibria $y_3=\pm 1$.
\begin{rem} \label{rem:SF_parameter_R_1}
    Note that here the variable $\sigma_3$ changes the geometry of the critical manifold 
      $\mathcal{S}_3$
    and $r_1$ is the singular perturbation parameter. In terms of the original system \eqref{equ:rob_reaction_2D} this means that - loosely speaking - $\eps_2$ changes the geometry and $\eps_1$ is the singular perturbation parameter. We will see that these roles will be switched when we study the dynamics in region $B_3$.
\end{rem}
The eigenvalue of the linearization of the layer problem  \eqref{eq:K_13_blow-up_K^3_r_1=0} is $\lambda=-2y_3-s_3(c-\sigma_3^2)$. At the line $s_3=0$, $y_3=1$ we obtain 
$$\lambda=-2<0$$
and we conclude that the line $s_3=0$, $y_3=1$ is part of the attracting branch $\mathcal{S}^a_3$ of the critical manifold, which extends regularly into $s_3>0$, since  $s_3$ acts as a  regular perturbation parameter in \eqref{eq:crit_mfld_K13^3}.

The reduced flow on $\mathcal{S}^a_3$ is given by 
\begin{equation} \label{eq:K_13^3_r=0_reduced}
    \begin{aligned}
    \Dot{s}_3&=\frac{1}{2}s_3^2y_3^2\\
    \Dot{\sigma}_3&=-\frac{1}{2}\sigma_3s_3y_3^2.
    \end{aligned}
\end{equation}
For $s_3=0$ the reduced flow along $\mathcal{S}^a_3$ is stationary, whereas for $s_3>0$ the variable $s_3$ 
increases and $\sigma_3$ decreases, see   Figure \ref{fig:K_13_singular dynamics}.

For the remaining analysis of $\mathcal{S}^a$ close to the top of the sphere, we change to the scaling chart $\mathcal{K}_{12}^3$ where the dynamics is governed by
\begin{equation} \label{eq:K_13_blow-up_K^2_dynamics}
    \begin{aligned}
    y_2'&= -z_2-y_2^2-y_2c-\sigma_2^2y_2z_2 -r_1y_2\\
    z_2'&=r_1y_2^2\\
    \sigma_2'&=0,
\end{aligned}
\end{equation}

which is again a slow-fast system with singular perturbation parameter $r_1$. Note that system \eqref{eq:K_13_blow-up_K^2_dynamics} has a line of equilibria (independent of $r_1$) at $y_2=z_2=0$ which corresponds to the genuine equilibrium $Q$.

The layer problem is given by
\begin{equation} \label{eq:K_13_blow-up_K^2_r_1=0}
    \begin{aligned}
    y_2'&= -z_2-y_2^2-y_2c-\sigma_2^2y_2z_2\\
    z_2'&=0\\
    \sigma_2'&=0,
\end{aligned}
\end{equation}

with critical manifold
\begin{equation} \label{eq:crit_mfld_K13^2}
    \mathcal{S}_2=\{(y_2,z_2,\sigma_2)^T \in \R^3:-z_2-y_2^2-y_2c-\sigma_2^2y_2z_2=0\}.
\end{equation}
For $\sigma_2=0$, i.e., on the sphere, the critical manifold has the simple form
$$z_2=-y_2(y_2+c)$$
with a fold point at $y_2=-\frac{c}{2}$ and non-vanishing eigenvalue $\lambda=-2y_2-c$. We conclude that $z_2=-y_2(y_2+c),\, y_2\geq 0$ is part of the attracting branch $\mathcal{S}^a_2$ of the critical manifold and extends regularly into $\sigma_2>0$ since $\sigma_2$ is a regular perturbation parameter in \eqref{eq:crit_mfld_K13^2}.
The critical manifold $\mathcal{S}^a_2$ is uniformly normally attracting for $y \geq 0$ and $\sigma_2 \geq 0$ small enough since for $\sigma_2=0$ the fold point at $y_2=-\frac{c}{2}$ is bounded away from the half space $y_2\geq 0$.
The reduced flow on $S_2^a$ is given by 
$$\Dot{z}_2=y_2^2$$
such that orbits along $S_2^a$ with $y_2(0)>0$ converge to the line of equilibria $y_2=z_2=0$ corresponding to $Q$ in a center-like manner, i.e. with an algebraic rate.
By Fenichel theory we conclude that there exists a two-dimensional attracting invariant slow manifold $\mathcal{S}^a_{r_1}$ for $0<r_1 \ll 1$ with slow flow converging to the reduced flow on $\mathcal{S}^a$ as $r_1 \to 0$.
Since no new equilibria occur for $r_1 >0$ all orbits of the slow flow converge to a point on the line of equilibria corresponding to $Q$.

\end{proof}

Based on Lemma \ref{lem:S^a_B_12} we can now construct an $\Tilde{\eps}_2$-family of singular orbits connecting the line of initial values corresponding to $O$ with the line of equilibria corresponding $Q$.
 
\begin{lem} \label{lem:singular_orbit_B12}
    There exists a family of singular orbits $\gamma_0^{\Tilde{\eps}_2},\, \Tilde{\eps}_2 \in [0,1/\sqrt{\beta_2}]$ of the blown-up extended system of \eqref{eq:Bar{B}_1_rescaled_K13} connecting the line of initial values with the line of equilibria corresponding to $O$ and $Q$, respectively.
\end{lem}
\begin{proof}
We define the fast fibers connecting the line of initial conditions 
$$O_3= (0,\frac{\Tilde{\eps}_2}{\sqrt{c}},\sqrt{c})^T, \; \Tilde{\eps}_2 \in [0,1/\sqrt{\beta_2}]$$
with 
$$(y_3,s_3,\sigma_3)^T=(1,\frac{\Tilde{\eps}_2}{\sqrt{c}},\sqrt{c})^T \in \mathcal{S}_3^a$$ as $\gamma_0^{\Tilde{\eps}_2,f}$.
The orbits under the reduced flow along the critical manifold $\mathcal{S}^a$ connecting 
$$(y_3,s_3,\sigma_3)^T=(1,\frac{\Tilde{\eps}_2}{\sqrt{c}},\sqrt{c})^T \in \mathcal{S}_3^a$$
with 
the line of equilibria
$$Q_2=(0,0,\Tilde{\eps}_2)^T \in \R^3:\Tilde{\eps}_2 \in [0,1/\sqrt{\beta_2}]\}$$
are defined as $\gamma_0^{\Tilde{\eps}_2,s}$

The $\Tilde{\eps}_2$-family of singular orbits is therefore given by $$\gamma_0^{\Tilde{\eps}_2}:=\gamma_0^{\Tilde{\eps}_2,f} \cup \gamma_0^{\Tilde{\eps}_2,s},$$
see Figure \ref{fig:K_13_singular dynamics}.
\end{proof}

In the following result we prove that the singular orbits $\gamma_0^{\Tilde{\eps}_2}$ perturb to smooth orbits connecting $O$ and $Q$ for $0<r_1 \ll 1$.

\begin{thm} \label{thm:K_12}
    There exists a constant $\Tilde{r}>0$ such that for all $\Tilde{\eps}_2 \in (0,1/\sqrt{\beta_2}]$ and $r_1 \in (0,\Tilde{r}]$, there exists a smooth orbit $\gamma_{r_1}^{\Tilde{\eps}_2}$ of \eqref{eq:Bar{B}_1_rescaled_K13} connecting the initial value $O$ with the genuine equilibrium $Q$. The corresponding orbit in blown-up space $\Bar{\gamma}_{r_1}^{\Tilde{\eps}_2}$ is $\mathcal{O}(r_1)$-close to its corresponding singular orbit $\gamma_0^{\Tilde{\eps}_2}$ in Hausdorff distance.
\end{thm}
\begin{proof}
The existence of the   singular orbits in Lemma \ref{lem:singular_orbit_B12}, 
    standard Fenichel theory,  Lemma \ref{lem:S^a_B_12} and arguments similar to the proof of
    Theorem \ref{thm:convergence_R_2} imply that the forward solution   with initial value $O$ 
   converges to $Q$ for $0 < r_1 \ll 1$.  We denote this solution by $\Bar{\gamma}_{r_1}^{\Tilde{\eps}_2}$ which by construction is
   $\mathcal{O}(r_1)$-close to  the singular orbit $\Bar{\gamma}_{0}^{\Tilde{\eps}_2}$
    for all $\Tilde{\eps}_2 \in (0,1/\sqrt{\beta_2}]$ and $0<r_1 \ll 1$. The assertions of the theorem follow by applying the blow-up transformation \eqref{eq:spherical_blow-up_R1}, i.e., $\gamma_{r_1}^{\Tilde{\eps}_2}=\Phi(\Bar{\gamma}_{r_1}^{\Tilde{\eps}_2})$.
\end{proof}

\begin{rem}
    Note that $\Tilde{\eps}_2=0$ is not included in Theorem \ref{thm:K_12}, since this corresponds to the original parameters $\eps_1=\eps_2=0$. In this case we do not observe dynamics because $y=0$ is a line of equilibria, as already mentioned in the rough classification in the beginning of Section \ref{sec:Parameter blow-up}.
\end{rem}

This concludes the analysis in region $\Bar{B}_1$. It remains to investigate the dynamics in region $\Bar{B}_3$.

\section{Analysis in Region $B_3$} \label{sec:Regime R3}
The analysis in region $\Bar{B}_3$ is carried out in chart $\mathcal{P}_2$. Inserting \eqref{equ:parameter_blow-up-chart_K3} into \eqref{equ:rob_reaction_2D} we obtain
\begin{equation} \label{eq:K_3}
    \begin{aligned}
        y'&=r^2\tilde{\eps}_1(c-y-z)-y^2-ryz\\
        z'&=y^2.
    \end{aligned}
\end{equation}
For $r=0$ this results in the same limiting system as in chart $\mathcal{P}_1$, see \eqref{eq:degenerate_limit} and Figure  \ref{fig:rob_reaction_2d_k1=k3=0}, with non-hyperbolic critical manifold 
$$y=0.$$
Rescaling  $y$ with \eqref{eq:rescaling_R2}, as before, we obtain (after dividing out a factor of $r$)
\begin{equation} \label{eq:K_3_rescaled}
    \begin{aligned}
        \Tilde{y}'&=\tilde{\eps}_1(c-r\Tilde{y}-z)-\Tilde{y}^2-\Tilde{y}z\\
        z'&=r\Tilde{y}^2.
    \end{aligned}
\end{equation}

System \eqref{eq:K_3_rescaled} is of standard slow-fast type with singular perturbation parameter $r$. In the following we will again omit the \enquote{ $\Tilde{ }$ }. 

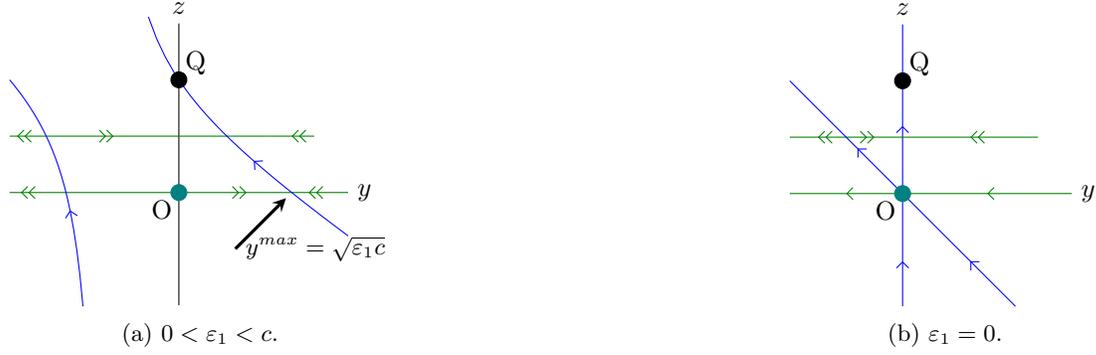
\begin{figure}
    \begin{subfigure}{0.4\textwidth}
        \centering
			\begin{tikzpicture}[scale=1.5]
			
\draw[] (0,-1) -- (0,1.5,0) node[above]{\textcolor{black}{$z$}};
        \draw[->>-=.8,black!50!green] (-1,0) -- (1,0);
        \draw[->>-=.8,black!50!green] (-1,0) -- (-1.5,0);
        \draw[->>-=.7,black!50!green] (1.5,0) node[right]{\textcolor{black}{$y$}} -- (1,0) ;
        \draw[->>-=.5,black!50!green] (1.2,0.5) -- (0.8,0.5);
        \draw[->>-=.3,black!50!green] (-1.17,0.5) -- (0.8,0.5);
        \draw[->>-=.8,black!50!green] (-1.17,0.5) -- (-1.5,0.5);
        \draw[->-=.4,scale=1, domain=1.5:-0.27, variable=\x, blue] plot ({\x}, {(1-\x*\x)/(1+1.5*\x)});
        \draw[->-=.4,scale=1, domain=-0.85:-1.5, variable=\x, blue] plot ({\x}, {(1-\x*\x)/(1+1.5*\x)});
        \filldraw[teal] (0,0) circle (2pt);
        \filldraw[black] (0,1) circle (2pt);
        \draw[] (-0.15,-0.15) node[] {O};
        \draw[] (0.15,1.15) node[] {Q};
        \draw[-stealth,very thick,black] (0.5,-0.5) node[right,black] {\small $y^{max}=\sqrt{\eps_1c}$} -- ++ (0.45,0.45) ;
        
        ;
			
			\end{tikzpicture}
			\caption{$0<\eps_1<c$.}
			\label{fig:K3_rescaled_singular_k>0}
    \end{subfigure}
    \hfill
    \begin{subfigure}{0.4\textwidth}
        \centering
			\begin{tikzpicture}[scale=1.5]
			
\draw[->-=.4,blue] (0,-1) -- (0,0) ;
\draw[->-=.4,blue] (0,0) -- (0,1.5) node[above]{\textcolor{black}{$z$}};
        \draw[->-=.5,black!50!green] (0,0) -- (-1,0);
        \draw[->-=.5,black!50!green] (1.5,0) node[right]{\textcolor{black}{$y$}} -- (0,0) ;
        \draw[->>-=.5,black!50!green] (1.2,0.5) -- (0,0.5);
        \draw[->>-=.5,black!50!green] (-0.5,0.5) -- (0,0.5);
        \draw[->>-=.5,black!50!green] (-0.5,0.5) -- (-1,0.5);
        \draw[->-=.4,scale=1, domain=1:0, variable=\x, blue] plot ({\x}, {-\x});
        \draw[->-=.4,scale=1, domain=0:-1, variable=\x, blue] plot ({\x}, {-\x});
        \filldraw[teal] (0,0) circle (2pt);
        \filldraw[black] (0,1) circle (2pt);
        \draw[] (-0.15,-0.15) node[] {O};
        \draw[] (0.15,1.15) node[] {Q};
        
        ;
			
			\end{tikzpicture}
			\caption{$\eps_1=0$.}
			\label{fig:K3_rescaled_singular_k=0}
    \end{subfigure}
    \caption{Singular dynamics of \eqref{eq:K_3_rescaled}.}
    \label{fig:K3_rescaled_singular}
\end{figure}

The corresponding layer problem is given by 
\begin{equation} \label{eq:K_3_rescaled_r=0}
    \begin{aligned}
        y'&=\eps_1(c-z)-y^2-yz\\
        z'&=0.
    \end{aligned}
\end{equation}

which for $\eps_1 >0$, resembles (actually is identical to) the situation in region $\Bar{B}_2$. Indeed, the critical manifold is given by 
\begin{equation} \label{K_3_rescaled_r=0_crit_mfld}
    \mathcal{S}=\{(y,z)^T \in \R^2: \eps_1(c-z)-y^2-yz=0\}
\end{equation}
and is normally attracting (repelling) for all $y>-\eps_1$ ($y<-\eps_1)$), see Figure \ref{fig:K3_rescaled_singular_k>0} and compare with Figure \ref{fig:K1_R2_rescaled_singular2}.

In the limit $\eps_1 \to 0$ normal hyperbolicity is lost at the origin since
for $\eps_1=0$ the critical manifold consists of two lines 
$$y=0 \text{ and } z=-y$$
which intersect at the origin. The linearization of the layer problem \ref{eq:K_3_rescaled_r=0}
at these lines has eigenvalue $\lambda_1=-z$ and $\lambda_2=-y$, respectively. Hence, the critical manifold $\mathcal{S}$ is not normally hyperbolic at the origin for $\eps_1=0$, see Figure \ref{fig:K3_rescaled_singular_k=0}.

To regain normal hyperbolicity we once again enlarge phase space by adding the  equation $\eps_1'=0$
and blow-up the degenerate equilibrium 
$(y,z,\eps_1)^T =(0,0,0)^T$
of this extended system. The suitable blow-up transformation is
\begin{equation} \label{eq:spherical_blow-up_R3}
     \begin{aligned}
         \Phi:[0,\infty) \times \mathbb{S}^2  & \to \R^3\\
	(\sigma,\bar{y},\bar{z},\bar{\eps}_1) & \mapsto \begin{cases}
	    y=\sigma \bar{y}\\
        z=\sigma\bar{z}\\
        \eps_1=\sigma^2 \bar{\eps}_1,
	\end{cases}
     \end{aligned}
 \end{equation}
which uses the same weights as the
analysis of the slow passage through a transcritical bifurcation, see \cite{krupa_extending_transcritical}.

\begin{lem} \label{S^a_B3}
    In the blown-up space of system \eqref{eq:K_3_rescaled} extended by the equation $\eps_1'=0$ there exists for $r=0$ a two-dimensional attracting critical manifold $\mathcal{S}^a$ (shown blue in Figure \ref{fig:K_3_singular dynamics}), which contains the line of equilibria corresponding to the genuine equilibrium $Q$. The critical manifold $\mathcal{S}^a$ perturbes regularly to a slow manifold $S_r^a$ for $r>0$ small enough. All orbits of the reduced flow on $\mathcal{S}^a$ approach the line of equilibria in the center direction, i.e., with an algebraic rate, for all $0<\eps_1\leq \beta_3$. 
\end{lem}

\begin{rem}
For better visibility we have changed the orientation in Figure \ref{fig:K_3_singular dynamics}, i.e. we look towards the origin from the $\bar{z} =1$ side of the sphere. 
\end{rem}

\begin{proof}
    Again it will be convenient to work in directional charts, which we denote by $\mathcal{K}_{3}^2$ and $\mathcal{K}_{3}^3$. The blow-up transformation in these charts is given by 
   \begin{align} 
        &\mathcal{K}_{3}^2: y=\sigma_2y_2, \quad z=\sigma_2z_2, \quad \eps_1=\sigma_2^2 \label{eq:spherical_blow-up_K_3^2}\\
        &\mathcal{K}_{3}^3: y=\sigma_3y_3, \quad z=\sigma_3, \quad \eps_1=\sigma_3^2\eps_{13},\label{eq:spherical_blow-up_K_3^3}
    \end{align}
covering  the top $\Bar{\eps}_1 >0$ and the front $\Bar{z} > 0$ part of the sphere, respectively,
see  Figure \ref{fig:K_3_singular dynamics}. The parts of $\mathcal{S}^a$ investigated in chart $\mathcal{K}_{3}^2$ and $\mathcal{K}_{3}^3$ are denoted by $\mathcal{S}^a_2$ and $\mathcal{S}^a_3$, respectively.
Since we have blown-up the initial value $O$, we start the analysis in the scaling chart $\mathcal{K}_3^2$, where the dynamics is governed by 
\begin{equation} \label{eq:K_3^2_dynamics}
    \begin{aligned}
        y_2'&=c-\sigma_2z_2-y_2^2-y_2z_2-r\sigma_2y_2\\
        z_2'&=ry_2^2\\
        \sigma_2'&=0.
    \end{aligned}
\end{equation}

System \eqref{eq:K_3^2_dynamics} is of standard slow-fast type with singular perturbation parameter $r$. The corresponding layer problem is given by 
\begin{equation} \label{eq:K_3^2_r=0}
    \begin{aligned}
        y_2'&=c-\sigma_2z_2-y_2^2-y_2z_2\\
        z_2'&=0\\
        \sigma_2'&=0
    \end{aligned}
\end{equation}
For $r=0$ we find the critical manifold
\begin{equation} \label{K_3^2_r=0_crit_mfld}
    \mathcal{S}_2=\{(y_2,z_2,\sigma_2)^T \in \R^3: c-\sigma_2z_2-y_2^2-y_2z_2=0 \},
\end{equation}
which simplifies on the sphere $\sigma_2=0$ to $z_2=\frac{c}{y_2}-y_2$.
   \begin{figure}
    \centering
    \includegraphics{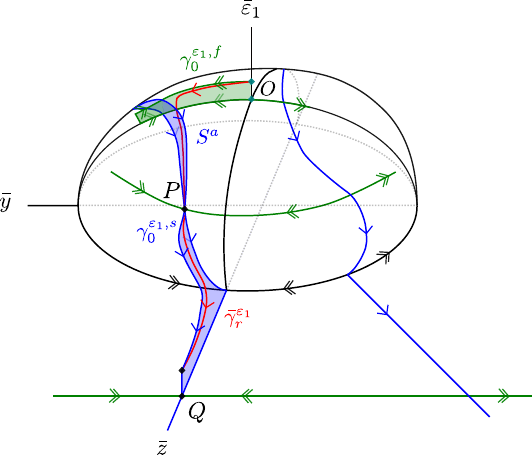}
    \caption{Dynamics of the blown-up extended system \eqref{eq:K_3_rescaled}.}
    \label{fig:K_3_singular dynamics}
\end{figure}

As already indicated in Remark \ref{rem:SF_parameter_R_1} we note the following
\begin{rem} \label{rem:SF_parameter_R_3}
    In system \eqref{eq:K_3^2_dynamics} $r$ is the slow-fast parameter and $\sigma_2$ changes the geometry of the critical manifold. Translated to the original parameters, i.e., undoing the blow-up transformations \eqref{equ:parameter_blow-up-chart_K3} and \eqref{eq:spherical_blow-up_K_3^2}, this implies that in region $B_3$ we have $\eps_2$ as slow-fast parameter and $\eps_1$ changes the geometry of the critical manifold.
\end{rem}

The eigenvalue of the linearization of the layer problem \eqref{eq:K_3^2_r=0} is $\lambda=-2y_2-z_2$. We conclude that the curve $z_2=\frac{c}{y_2}-y_2,\, y_2>0$ on the sphere $\sigma_2=0$ is part of the normally attracting branch of the critical manifold $\mathcal{S}^a_2$, which again extends regularly into $\sigma_2>0$, since $\sigma_2$ acts as a regular perturbation parameter in \eqref{K_3^2_r=0_crit_mfld}. 

The reduced flow on $S_2^a$ is given by
    $$\Dot{z}_2=y_2^2,$$
hence  $z_2$ increases for all $y_2\neq 0$, see Figure \ref{fig:K_3_singular dynamics}.

For the remaining analysis of $\mathcal{S}^a$ away from the sphere, we change to the exit chart $\mathcal{K}_3^3$ where the dynamics is governed by
\begin{equation} \label{eq:K_3^3_dynamics}
    \begin{aligned}
        y_3'&=\eps_{13}(c-r\sigma_3y_3-\sigma_3)-y_3^2-y_3-ry_3^3\\
        \sigma_3'&=r\sigma_3y_3^2\\
        \eps_{13}'&=-2r\eps_{13}y_3^2.
    \end{aligned}
\end{equation}
The layer problem is now given by 
\begin{equation} \label{eq:K_3^3_r=0}
    \begin{aligned}
        y_3'&=\eps_{13}(c-\sigma_3)-y_3^2-y_3\\
        \sigma_3'&=0\\
        \eps_{13}'&=0,
    \end{aligned}
\end{equation}
with critical manifold
\begin{equation} \label{K_3^3_r=0_crit_mfld}
    \mathcal{S}_3=\{(y_3,\sigma_3,\eps_{13})^T \in \R^3: \eps_{13}(c-\sigma_3)-y_3^2-y_3=0\}.
\end{equation}

On the invariant plane $\eps_{13}=0$, the critical manifold $\mathcal{S}_3$ corresponds to the lines $y_3=0$ and $y_3=-1$. The eigenvalue of the linearization of the layer problem \eqref{eq:K_3^3_r=0} is $\lambda=-2y_3-1$, hence the line $y_3=0, \, \eps_{13}=0$ is part of the normally attracting branch $S_3^a$ of the critical manifold. As before this line extends regularly into $\eps_{13}>0$ to a smooth manifold $\mathcal{S}^a_3$.

The reduced flow on $S_3^a$ is given by 
\begin{equation} \label{eq:K_3^3_r=0_reduced}
    \begin{aligned}
        \sigma_3'&=\sigma_3y_3^2\\
        \eps_{13}'&=-2\eps_{13}y_3^2,
    \end{aligned}
\end{equation}
such that $\sigma_3$ increases and $\eps_{13}$ decreases for $y_3>0$ on $\mathcal{S}^a_3$. All orbits on $\mathcal{S}^a_3$ with $\sigma_3,\,\eps_{13}>0$ approach the line of equilibria $y_3=0$, $\sigma_3=\sqrt{c}$ corresponding to $Q$ (which is contained in $\mathcal{S}^a_3$) in a center-like manner, i.e., with an algebraic rate. For completeness note that on the invariant plane $\eps_{13}=0$ the critical manifold corresponds to the line $y_3=0$ with stationary reduced flow.
We conclude that there exists a two-dimensional attracting invariant slow manifold $\mathcal{S}_r^a$ for $0<r \ll 1$ with slow flow converging to the reduced flow on $\mathcal{S}^a$ as $r \to 0$. Since no new equilibria occur for $r>0$, all orbits of the slow flow converge to a point on the line of equilibria corresponding to $Q$. 

\end{proof}

Based on Lemma \ref{S^a_B3} we can now construct an $\eps_1$-family of singular orbits connecting the line of initial values corresponding to $O$ with the line of equilibria corresponding to $Q$.

\begin{lem} \label{lem:singular_orbit_B3}
    There exists a family of singular orbits $\gamma_0^{\eps_1},\, \eps_1 \in [0,\beta_3]$ of the blown-up extended system of \eqref{eq:K_3_rescaled} connecting the line of initial values with the line of equilibria corresponding to $O$ and $Q$, respectively.
\end{lem}
\begin{proof}
We define the fast fibers connecting the line of initial conditions 
$$O_2= (0,0,\sqrt{\eps_1})^T, \; \eps_1 \in [0,\beta_3]$$
with 
$$(y_2,z_2,\sigma_2)^T=(\sqrt{c},0,\sqrt{\eps_1})^T \in \mathcal{S}_2^a$$ as $\gamma_0^{\eps_1,f}$.
The forward orbits under the reduced flow along the critical manifold $\mathcal{S}^a$ connecting 
$$(y_2,z_2,\sigma_2)^T=(\sqrt{c},0,\sqrt{\eps_1})^T \in \mathcal{S}_2^a$$
with 
the line of equilibria
$$Q_3=\Big(0,c,\frac{\eps_1}{c^2}\Big)^T \in \mathcal{S}_3^a, \, \eps_1 \in [0,\beta_3]$$
are denoted as $\gamma_0^{\eps_1,s}$.
The $\eps_1$-family of singular orbits is therefore given by $$\gamma_0^{\eps_1}:=\gamma_0^{\eps_1,f} \cup \gamma_0^{\eps_1,s},$$
see Figure \ref{fig:K_3_singular dynamics}.
\end{proof}

The following theorem assures that the singular orbits $\gamma_0^{\eps_1}$ perturb to smooth orbits connecting $O$ and $Q$ for $0<r \ll 1$.

\begin{thm} \label{thm:K_3}
    There exists a constant $r_0>0$ such that for all $\eps_1 \in (0,\beta_3]$ and $r \in (0,r_0]$, there exists a smooth orbit $\gamma_r^{\eps_1}$ of \eqref{eq:K_3_rescaled} connecting the initial value $O$ with the genuine equilibrium $Q$. The corresponding orbit in blown-up space $\Bar{\gamma}_r^{\eps_1}$ is $\mathcal{O}(r)$-close to its corresponding singular orbit $\gamma_0^{\eps_1}$ in Hausdorff distance.
\end{thm}
\begin{proof}
    Combining Lemma \ref{S^a_B3} and Lemma \ref{lem:singular_orbit_B3} it follows from standard Fenichel theory with slow-fast parameter $r$ applied to the blown-up extended system of \eqref{eq:K_3_rescaled}  and arguments similar to the proof of Theorem \ref{thm:convergence_R_2} imply 
    that  the singular orbits $\gamma_0^{\eps_1}$ perturb to smooth orbits $\Bar{\gamma}_r^{\eps_1}$ converging to $Q$,
    which are $\mathcal{O}(r)$-close for all $\eps_1 \in (0,\eps_3]$ and $r>0$ small enough. The assertions of the theorem follow by applying the blow-up map \eqref{eq:spherical_blow-up_R3}, i.e.,  $\gamma_r^{\eps_1}=\Phi(\Bar{\gamma}_r^{\eps_1})$.
\end{proof}

\begin{rem}
    Note that $\eps_1=0$ (actually $\Tilde{\eps}_1=0$ since the theorem is stated in chart $\mathcal{P}_2$) is not included in Theorem \ref{thm:K_3}, since this corresponds to the original parameter $\eps_1=0$. In this case we do not observe dynamics because $y=0$ is a line of equilibria, as already mentioned in the rough classification in the beginning of Section \ref{sec:Parameter blow-up}.
\end{rem}

We conclude with the proof of the main result, i.e., Theorem \ref{thm:main_summary}.

\begin{proof}[Proof of Theorem \ref{thm:main_summary}]
    It follows from Theorems \ref{thm:convergence_R_2}, \ref{thm:P11}, \ref{thm:K_12}, \ref{thm:K_3} that in each of the regions $B_{11}$, $B_{12}$, $B_2$, $B_3$ there exists a different slow-fast structure of \eqref{equ:rob_reaction_2D} with a corresponding singular orbit $\gamma_0$ which perturbs to a genuine orbit for $\eps_1, \eps_2$ small. The error estimates in $\eps_1$ and $\eps_2$ for each case are obtained by undoing the rescalings of the blow-up transformations as in \eqref{eq:Evaluate_numerics_K_2}.
\end{proof}

\section{Summary and Outlook} \label{sec:Summary}
\begin{figure}[]
    \centering
    \includegraphics[width=0.7\textwidth]{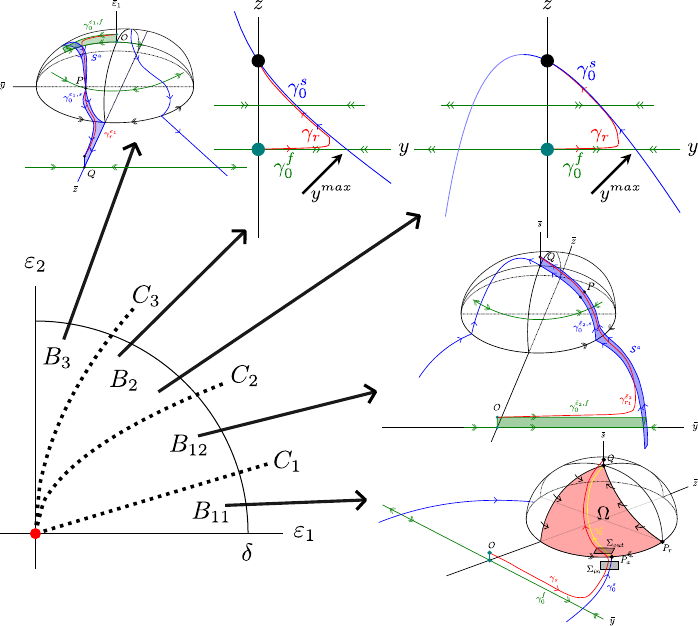}
    \caption{Multi-scale structure and singular orbits of \eqref{equ:rob_reaction_2D} in different regions of parameter space.}
    \label{fig:phase_plane_comparison}
\end{figure}

In this article, we conducted an asymptotic analysis of the Robertson model, a prominent example of stiffness in ODEs characterized by three reaction rates $k_1$, $k_2$, and $k_3$ of widely differing orders of magnitude. We focused on the scenario where $k_1,\, k_3 \ll k_2$. By rescaling the problem in terms of the small parameters $(\eps_1,\eps_2):=(k_1/k_2,k_3/k_2)$, we transformed the original equations into a two-parameter singular perturbation problem. To deal with the singular structures associated with the two small parameters, we introduced suitable blow-up transformations in parameter space. This allowed us to systematically explore the behaviour of the system in a neighbourhood of the singular limit $(\eps_1,\,\eps_2)=(0,0)$. Our analysis revealed four distinct scaling regimes  with different singular structures. Within each regime, we applied GSPT and further blow-ups in phase space to investigate the dynamics and the structure of the solutions. This combined approach enabled us to capture the various multi-scale structures of the model, see Figure \ref{fig:phase_plane_comparison}, which illustrates the main result Theorem \ref{thm:main_summary} and the details of the analysis given in Sections \ref{sec:Regime R2}, \ref{sec:Regime R1}, and \ref{sec:Regime R3}.
In each region we identified a specific type of singular orbit connecting the initial value $O$ with equilibrium $Q$, which perturbs to a genuine orbit (shown in red) for $\eps_1, \eps_2$ small.

\begin{table}[h!]
\setlength{\tabcolsep}{1em} 
{\renewcommand{\arraystretch}{1.4} 
\begin{center}
\begin{tabular}{ |c||c|c| } 
 \hline
 \textbf{Region} & $\mathbf{y^{max}_{num}}$ & $\mathbf{y^{max}}$ \\ 
 \hhline{|=|=|=|}
 $B_{11}$& $2.16\cdot 10^{-2}$ & $2.2\cdot 10^{-2}+\mathcal{O}(10^{-4})$ \\ 
 \hline
 $B_{12}$ & $6.98\cdot 10^{-3}$ & $7.0\cdot 10^{-3}+\mathcal{O}(10^{-5})$ \\ 
  \hline
 $B_2$ (lower) & $7.05\cdot 10^{-4}$ & $7.07\cdot 10^{-4}+\mathcal{O}(10^{-7})$ \\
  \hline
   $B_2$ (upper) & $7.07\cdot 10^{-5}$ & $7.071\cdot 10^{-5}+\mathcal{O}(10^{-9})$\\
  \hline
 $B_3$ & $2.26\cdot 10^{-6}$ & $2.236\cdot 10^{-6}+\mathcal{O}(10^{-10})$\\
  \hline
\end{tabular}
\end{center}}
    \caption{\comm{Comparison of the numerical ($y^{max}_{num}$) and the analytical ($y^{max}$) plateau value obtained at representative points $(\eps_1,\eps_2)$ in the different parameter regions corresponding to the simulations shown in Figure \ref{fig:numerical_comparison}.}}
    \label{tab:Value_comparison}
\end{table}

The asymptotic results derived from our analysis are in excellent qualitative and quantitative agreement with numerical simulations, compare Figure \ref{fig:numerical_comparison} with Figure \ref{fig:phase_plane_comparison}, e.g., the maximal value of the $y$-component. Our analysis predicts $y^{max}=\sqrt{\eps_1c}+\mathcal{O}(\eps_1)$ in $B_{11}$, $B_{12}$, and $B_{2}$ and $y^{max}=\sqrt{\eps_1c}+\mathcal{O}(\sqrt{\eps_1} \eps_2)$ in $B_3$. \comm{The values obtained for $y^{max}$ (with the same parameter values as in the simulations shown in Figure \ref{fig:numerical_comparison}) are collected in Table \ref{tab:Value_comparison}, which fit well with the corresponding numerical results $y_{num}^{max}$ of the simulations shown in Figure \ref{fig:numerical_comparison}.} 
Additionally we observe in Figure \ref{fig:numerical_comparison} that the time it takes for  $z$ to increase becomes longer as we move counterclockwise. \comm{The analytical explanation for this observation is the change in the geometry of the attracting part of the critical manifold $\mathcal{S}^a$ in the different regions, cf., Figure \ref{fig:phase_plane_comparison}.}



Overall, this work provides a thorough understanding of the dynamics and detailed asymptotics of the Robertson model. This case study  highlights the potential of combining GSPT with blow-up in parameter space for analyzing multi-parameter singular perturbation problems. We believe that this approach is applicable to more complicated problems and has the potential to lead to a framework for the analysis of multi-parameter 
singular perturbations. \comm{In ongoing work we are using this approach in the analysis of a five-variable model of the cell cycle \cite{deso,tyson1991} with singular dependence on three parameters.}

\clearpage

\section*{Acknowledgements} \label{sec:Ack}
\comm{We thank two anonymous referees for their helpful suggestions. The authors also acknowledge TU Wien Bibliothek for financial support through its Open Access Funding Programme.}

\printbibliography[heading=bibintoc,
title={References}]

\end{document}